\documentclass[a4paper, 12pt]{article}

\usepackage[sort&compress]{natbib}
\bibpunct{(}{)}{;}{a}{}{,} 

\usepackage{amsthm, amsmath, amssymb, mathrsfs, multirow, url, subfigure}
\usepackage{graphicx} 
\usepackage{ifthen} 
\usepackage{amsfonts}
\usepackage[usenames]{color}
\usepackage{fullpage}
\usepackage{tikz}
\usetikzlibrary{shapes.symbols, hobby, positioning}
\usepackage[ruled,vlined]{algorithm2e}

\theoremstyle{plain} 
\newtheorem{thm}{Theorem}
\newtheorem{cor}{Corollary}

\newtheorem*{prop0}{Proposition}

\theoremstyle{definition}

\theoremstyle{remark}

\newtheorem{remark}{Remark}

\newtheorem*{gammaex}{Gamma Example}
\newtheorem*{gaussex}{Normal Example}
\newtheorem*{question}{Question}
\newtheorem{trick}{Trick}

\newcommand{\prob}{\mathsf{P}} 
\newcommand{\E}{\mathsf{E}}

\newcommand{\bin}{{\sf Bin}}
\newcommand{\unif}{{\sf Unif}}
\newcommand{\nm}{{\sf N}}

\newcommand{\chisq}{{\sf ChiSq}}

\newcommand{\mult}{{\sf Mult}}

\newcommand{\RR}{\mathbb{R}}

\newcommand{\ZZ}{\mathbb{Z}}

\newcommand{\TT}{\mathbb{T}}

\renewcommand{\S}{\mathcal{S}}

\newcommand{\iid}{\overset{\text{\tiny iid}}{\,\sim\,}}
\newcommand{\ind}{\overset{\text{\tiny ind}}{\,\sim\,}}

\newcommand{\model}{\mathscr{P}}
\newcommand{\prior}{\mathsf{Q}}

\newcommand{\kernel}{\mathsf{K}}
\newcommand{\marg}{\mathsf{M}}

\newcommand{\cred}{\mathscr{C}}

\newcommand{\lPi}{\amalg}%{\underline{\Pi}}
\newcommand{\uPi}{\Pi}%{\overline{\Pi}}

\newcommand{\lprob}{\underline{\mathsf{P}}}
\newcommand{\uprob}{\overline{\mathsf{P}}}

\title{Possibilistic inferential models: a review}
\author{Ryan Martin\footnote{Department of Statistics, North Carolina State University, {\tt rgmarti3@ncsu.edu}}
}
\date{\today}

\begin{document}

\maketitle 

\begin{abstract}
An {\em inferential model} (IM) is a model describing the construction of provably reliable, data-driven uncertainty quantification and inference about relevant unknowns.  IMs and Fisher's fiducial argument have similar objectives, but a fundamental distinction between the two is that the former doesn't require that uncertainty quantification be probabilistic, offering greater flexibility and allowing for a proof of its reliability.  Important recent developments have been made thanks in part to newfound connections with the imprecise probability literature, in particular, possibility theory.  The brand of {\em possibilistic IMs} studied here are straightforward to construct, have very strong frequentist-like reliability properties, and offer fully conditional, Bayesian-like (imprecise) probabilistic reasoning.  This paper reviews these key recent developments, describing the new theory, methods, and computational tools.  A generalization of the basic possibilistic IM is also presented, making new and unexpected connections with ideas in modern statistics and machine learning, e.g., bootstrap and conformal prediction.   

\smallskip

\end{abstract}

%\tableofcontents 

\section{Introduction}
\label{S:intro}

An {\em inferential model} (IM) is a model for data-driven uncertainty quantification and inductive inference about relevant unknowns.  These unknowns might be parameters in a posited statistical model, or functions thereof, but other situations are possible; see Section~\ref{S:beyond}.  More specifically, an IM offers mathematically rigorous, fully conditional, Bayesian-like uncertainty quantification---without requiring prior distributions or Bayes's theorem---that is provably reliable in the frequentist-like sense that its output is naturally calibrated across repeated sampling.  At a high level at least, this is reminiscent of what Fisher's fiducial argument aimed to do, so it's important to make clear IM's novelty at the outset: there are limits to the reliability of probabilistic uncertainty quantification in the absence genuine prior information, and the IM framework sidesteps these limitations not by relaxing ``reliability'' but instead by relaxing ``probabilistic,'' by working in a more flexible (but still mathematically rigorous) framework that doesn't require single or precise probability values to be assigned to each proposition about the unknowns.  These details will be explored in-depth in the sections that follow.  For now, suffice it to say that this novel use of concepts and tools coming the theory of {\em imprecise probability} in the context of statistical inference is surely what inspired \citet{cui.hannig.im} to describe IMs as ``one of the original statistical innovations of the 2010s.''  

It's been about 10 years since the monograph {\em Inferential Models: Reasoning with Uncertainty} \citep{imbook} was published, and a lot of exciting developments have been made since then on all fronts: foundations, theory, methodology, and computation.  An impetus behind these more recent advancements was the recognition that, while IMs' deviation from ordinary probability theory puts one in unfamiliar territory, that new territory isn't uninhabited---there's an entire community of researchers and a body of literature devoted to imprecise probability theory ripe with important insights, understanding, and mathematical/computational tools.  These efforts have inspired other developments \citep[e.g.,][]{xie.wang.repro, williams.cpfid, caprio.etal.isipta25} and have created opportunities for new-and-improved understanding of Bayes and fiducial inference, bootstrap, conformal prediction, etc.   The goal of this paper to review these recent developments, which are related to but also different from those in the aforementioned monograph, and to present some new insights, methods, and results.  I hope that this survey will make these exciting developments more accessible and draw a new generation of researchers' attention to these foundational advancements and open problems.  

Towards this goal, the remainder of the paper is organized as follows.  Section~\ref{S:background} sets the scene by reviewing probabilistic uncertainty quantification and highlighting its shortcomings.  In particular, Section~\ref{SS:false.conf} offers a new perspective on the {\em false confidence theorem} \citep{balch.martin.ferson.2017} establishing that any data-dependent posterior probability distribution---Bayes with any prior, (generalized) fiducial, etc.---will, in vacuous prior information cases like I'm considering here, tend to assign high probability or {\em confidence} to some false hypotheses.  The implications of this are two-fold: that probabilistic uncertainty quantification has an inherent unreliability and that correcting this requires quantification of uncertainty using more flexible tools from the world of imprecise probability.  Section~\ref{S:details} dives into the details of the {\em possibilistic} IM construction, where I'm emphasizing ``possibilistic'' because this new perspective---which is similar to that presented in \citet{imbook}---leans heavily on possibility-theoretic interpretations, calculus, computational tools, and mathematical structure.  The IM's key properties are presented next, including a finite-sample validity property which ensures that its imprecise-probabilistic output is properly calibrated and, in particular, that tests and confidence sets derived from the IM's output have frequentist error rate guarantees.  Efficiency is also discussed, including a new possibilistic {\em Bernstein--von Mises theorem} ensuring that the IM's output is asymptotically possibilistic-Gaussian and efficient like in the classical sense where the asymptotic variance matches the Cram\'er--Rao lower bound.  But IMs are not solely frequentist---they also offer fully conditional quantification of uncertainty and, in Section~\ref{SS:conditional}, I survey the relevant Bayesian-like properties that have not received much attention in previous works.  IM computation is nontrivial, and Section~\ref{SS:computation} reviews briefly an exciting new development facilitating sampling-based Monte Carlo computations for evaluating the IM's non-probabilistic output.  Section~\ref{S:reimagined} explains how IMs offer much more than a ``unification'' of frequentism and Bayesianism---the IM framework actually fills the holes in both paradigms! Section~\ref{S:marginal} deals with the important practical problem of eliminating nuisance parameters, and this motivates some new developments presented in Section~\ref{S:beyond}, which help IMs reach the next level of uncertainty quantification beyond those cases involving a statistical model.  Applications of this idea to inference on risk minimizers and to (conformal) prediction are considered there.  The paper ends in Section~\ref{S:discuss} with a brief summary, mentioning some relevant topics that are not discussed in this review, along with some open questions for future investigation.  The appendix/supplement offers supporting technical details and some additional examples.  

Some may argue that foundational efforts like these are impractical, but I disagree.  For years we've heard that statistics as a discipline is at risk of ``missing the boat'' when it comes to data science; see, e.g., the recent report by \citet{he.etal.crossroads}.  Statisticians collectively agree that statistics is an important part of data science, so why are we so concerned?  A boat's captain isn't concerned about missing the boat, neither are the first and second mates; only ancillary crew members and passengers worry about missing the boat.  This fear of missing the data science boat betrays our community's deep-down insecurity---a fear that we're only ancillary crew members---which surely can't be a consequence of us not contributing to enough applied projects or not proving enough consistency theorems.  It must instead be due to a deeper, more significant shortcoming, say: ``How can a discipline, central to science and to critical thinking, have two methodologies, two logics, two approaches that frequently give substantially different answers to the same problems'' \citep{fraser2011.rejoinder}. 
%... and yet no sense of urgency to clarify the conflicts.'' 
Our discipline can't be confident about its contributions when such fundamental questions remains unanswered.  Orienting our discipline by settling these foundational questions would prove that we have unique expertise to contribute, thereby giving us confidence that the boat can't leave us behind.

\section{Background and motivation}
\label{S:background}

% Expanded version of this section makes a chapter in my textbook... 

\subsection{Problem setup and notation}

Let $Z$ denote the observable data, taking values in a sample space $\ZZ$.  It'll often be the case that $Z$ has components $Z_1,\ldots,Z_n$, where $Z_i = (X_i, Y_i)$ might be an independent--dependent variable pair, etc.  When it's relevant that the data is a sample of size $n$, e.g., when considering large-sample properties, I'll write $Z^n = (Z_1,\ldots,Z_n)$.  Throughout, $z$ (or $z^n$) will denote a particular realization of the observable data $Z$ (or $Z^n$).  

Next, suppose statistical model $\{\prob_\theta: \theta \in \TT\}$, consisting of probability distributions supported on (subsets of) $\ZZ$, is imposed to quantify the variability or aleatory uncertainty in the observable data $Z$; I'll relax this parametric model restriction in Section~\ref{S:beyond} and later.  The probability distributions $\prob_\theta$ have a corresponding density/mass function $p_\theta(z)$ and, if $Z=z$ is the observed data, then $\theta \mapsto L_z(\theta) := p_\theta(z)$ is the likelihood function.  

Finally, let $\Theta$ denote the uncertain true parameter value, i.e., the one such that the statement ``$Z \sim \prob_\Theta$'' is true; this too can and will be relaxed later.  Like with my $z$ and $Z$ notation, note that I'm writing $\theta$ for generic parameter values and $\Theta$ for the uncertain true value to be inferred.  Importantly, here I'll be assuming that prior information about $\Theta$ is {\em vacuous}, i.e., my answer to the question ``What is the prior probability that $\Theta$ is in $A$?''~would be ``between 0 and 1'' for all $A \not\in \{\varnothing, \TT\}$---complete ignorance.  While ignorance is the state most commonly assumed in the statistics literature, it's not because complete ignorance is realistic---arguably it's rare in applications to know {\em literally nothing} about the quantity of interest.  See Section~\ref{S:discuss} for a brief discussion about recent developments allowing for non-vacuous, incomplete prior information about $\Theta$.    

Since prior information is vacuous, all that's available is the model/likelihood for $Z$ and the realization $z$.  According to \citet{hacking.logic.book}, ``Statisticians want numerical measures of the degree to which data support hypotheses,'' which, to me, sounds like {\em probabilistic uncertainty quantification} (Section~\ref{SS:prob.uq}).  So, despite having fully vacuous prior information, which puts proper Bayesian inference out of reach, the goal still is to assign data-dependent probabilities (or something similar) to hypotheses about the unknown $\Theta$.  Towards this, I'll follow Fisher---``the world’s master of quantifying uncertainty'' \citep{pearl.why}---, Jeffreys, Dempster, Berger, Walley, and other leaders.

\subsection{Probabilistic uncertainty quantification}
\label{SS:prob.uq}

For now, probabilistic uncertainty quantification means assigning data-dependent probabilities, say, $\prior_z(\cdot)$ to hypotheses about the unknown $\Theta$.  The point is that a hypothesis ``$\Theta \in H$'' determined by a subset $H \subseteq \TT$ can be true or not, and a large $\prior_z(H)$ assigned to that hypothesis would naturally be interpreted as an indication that data $z$ supports the truthfulness of that hypothesis.  (Henceforth, when I refer to a ``hypothesis $H$'' I mean the hypothesis ``$\Theta \in H$.'')  So, the magnitudes of these probabilities can be broadly used to assess which hypotheses the data does and doesn't support.  Of course, these probabilities can be used for other purposes (e.g., credible/confidence sets, prediction, and decision-making) but it's this basic function of assessing broadly where the data lends its support that distinguishes probabilistic uncertainty quantification from other approaches.  In Section~\ref{SS:existing} below I give a brief summary of the two most familiar of the existing approaches to probabilistic uncertainty quantification.  

It's important to distinguish probabilistic uncertainty quantification about an unknown and unobservable $\Theta$, like what's in consideration here, from that about an unknown that's observable, e.g., a future data point.  In the latter case, the probability model in question can be directly tested against observations: if a pre-specified event that the model says has (effectively) zero probability happens, then the model must be wrong.  This is {\em Cournot's principle}, see \citet{vovk1993}, \citet{shafer2007}, and \citet[][Ch.~10]{shafer.vovk.book.2019}.  In the former case, however, the true $\Theta$ will typically never be revealed, so the probabilistic uncertainty quantification about $\Theta$ can't be {\em directly} tested against reality.  But it can be {\em indirectly} scrutinized for reliability:
%Motivation for my perspective here comes from \citet{reid.cox.2014}:
\begin{quote}
Even if an empirical frequency-based view of probability is not used directly as a basis for inference, it is unacceptable if a procedure...~of representing uncertain knowledge would, if used repeatedly, give systematically misleading conclusions. \citep{reid.cox.2014} 
\end{quote}
In the present context, the ``procedure'' in question is a mapping $z \mapsto \prior_z$ from data $z$ to a probability distribution $\prior_z$ on $\TT$, and such procedure gives systematically misleading conclusions if $\prior_Z(H)$ tends to be large, as a function of $Z \sim \prob_\Theta$, when the hypothesis is false, i.e., $H \not\ni \Theta$.  So, at a high level, a probabilistic quantification of uncertainty is reliable if the ``objective'' $\prob_\Theta$-probability that $\prior_Z(H)$ is large for false hypotheses $H$ is itself small.  In symbols, but still only roughly, if I define the function 
\[ H \mapsto \sup_{\theta \not\in H} \prob_\theta\{ \text{$\prior_Z(H)$ is large} \}, \quad H \subset \TT, \]
then it's clearly a desirable property of $z \mapsto \prior_z$ that the right-hand side above be small for all $H \subset \TT$.  Importantly, a more precisely defined version of this function (see Section~\ref{SS:false.conf}) can be evaluated and used to assess the reliability of a given probabilistic uncertainty quantification procedure.  That these considerations assess the procedure's {\em reliability} is implied by Cournot's principle, albeit from a different perspective than that above.  Indeed, since events with small $\prob_\Theta$-probability effectively don't happen, if the right-hand side of the above display is small in the sense above, then $\prior_Z$ won't assign high probability to false hypotheses, hence no ``systematically misleading conclusions.''

\subsection{Existing approaches}
\label{SS:existing}

\paragraph{Default-prior Bayes.}
The idea of using equal probabilities as a default when genuine information is lacking has a long history.  It appears in the original work of \citet{bayes1763} and was adopted by \citet{laplace} and other contemporaries, eventually acquiring the name {\em principle of insufficient reason} \citep[e.g.,][p.~127--129]{stigler1986}.  \citet[][Ch.~4]{keynes.probability} later renamed it the {\em principle of indifference}, which he described as follows:
\begin{quote}
The Principle of Indifference asserts that...
%if there is now {\em known} reason for predicating of our subject one rather than another of several alternatives, then relatively to such knowledge the assertions of each of these alternatives have an {\em equal} probability.  Thus 
{\em equal} probabilities must be assigned to each of several arguments, if there is an absence of positive ground for assigning {\em unequal} ones. \citep[][p.~45]{keynes.probability}
\end{quote}
On the one hand, the above principle seems generally agreeable at least at first consideration, and it has been used in all sorts of applications and generalized in various ways \citep[e.g.,][]{jaynes2003}.  On the other hand, the principle of indifference has been heavily criticized by many authors, including Keynes and Fisher.  

Jeffreys responded to Fisher's criticism with a different perspective.  Rather than attempting to give a probabilistic description of ignorance, an impossible task (see Section~\ref{SS:false.conf}), he focused on the construction of otherwise justifiable default priors:
\begin{quote}
... find a way of saying that the magnitude of a parameter is unknown, when none of the possible values need special attention \citep[][p.~117]{jeffreys1961}.
\end{quote}
His efforts resulted in the now widely used Jeffreys priors \citep{jeffreys1946}, later shown to yield posterior distributions with excellent large-sample properties \citep[e.g.,][]{welch.peers.1963, datta.ghosh.1995}.  Substantial generalizations of Jeffreys formulation, in various directions, are now available \citep[e.g.,][]{berger.objective.book}. But despite this progress, there's still no general consensus about which---if any---default priors are ``right,'' so apparently the fundamental problem remains unsolved; see Efron's quote on page~\pageref{quote:efron}.

\paragraph{Fiducial and the like.}
A novel, non-Bayesian approach to probabilistic uncertainty quantification was taken by \citet{fisher1930, fisher1933, fisher1935a, fisher1935b}---what \citet{savage1961} famously described as ``a bold attempt to make the Bayesian omelet without breaking the Bayesian egg.'' I won't get into details here, so I refer the reader to \citet{zabell1992} and \citet{savage1976} for more on Fisher's ideas and to \citet{xie.singh.2012}, and \citet{hannig.review}, \citet{schweder.hjort.book} for some modern perspectives on fiducial-like inference. 

Roughly, Fisher's fiducial argument takes the model-based, parameter-dependent probabilities assigned to events about the observable data, reinterprets the events as (data-dependent) assertions about the unknown parameter, and then flips the aforementioned probabilities assigned to these events into subjective probabilities about the unknown parameter, given the observed data.  
The adjective ``fiducial''---which means founded on faith or trust---that Fisher chose to describe his solution makes it clear that he understood his was not a 100\% sound mathematical argument.  Fisher must have had in mind some principle that justifies the aforementioned faith/trust in his fiducial probabilities but, to my knowledge, he didn't clearly state any such principle.  \citet{dempster1963, dempster1964} described it as a {\em continue to regard} operation;
%, i.e., treat $F_\theta(z)$ at the given $Z=z$ and unknown $\theta$ exactly as one treated $F_\theta(Z)$ as a function of $Z \sim \prob_\theta$ for known $\theta$; 
Hannig and others \citep[e.g.,][]{murph.etal.fiducial, hannig.review} coined it a {\em switching principle} where what's random and what's fixed is switched.
%; see, also, Remark~\ref{re:interpretation} in Appendix~\ref{A:remarks}.
%and defining a subjective probability about $\Theta$, given $Z=z$, to match the ``objective'' probability about $Z$, given $\Theta$, bears some resemblance to the {\em principal principle} in \citet{lewis1980}.  
In any case, the fiducial argument involves a mix of mathematical reasoning and application of principles, hence is not so much different than default-prior Bayes.  

Thanks to Fisher's fame and the mystery surrounding his proposed solution, the fiducial argument received a lot of attention---and heavy scrutiny.  Damaging blows to Fisher's formulation were delivered by \citet{lindley1958}, \citet{dempster1963, dempster1964}, and \citet{buehler.fedderson.1963}, but these demonstrations, while insightful, largely only confirm that the fiducial argument is not mathematically sound.  
%Authors often quote \citet{efron1998} in calling fiducial inference ``Fisher's biggest blunder,'' but this is too negative and, out of context, doesn't accurately reflect Efron's point. 
%The goal of research is to advance the field, and Fisher's efforts did just that; for example, confidence limits are a simple reformulation of Fisher's proposal \citep{neyman1934, neyman1941}. 
That the solution offered by Fisher lacks mathematical rigor doesn't mean that the problem is impractical, immaterial, or impossible, hence fiducial inference remains a sort of Holy Grail for statisticians:
\begin{quote} \label{quote:zabell}
Fisher's attempt to steer a path between the Scylla of unconditional behaviorist methods which disavow any attempt at ``inference'' and the Charybdis of subjectivism in science was founded on important concerns, and his personal failure to arrive at a satisfactory solution to the problem means only that the problem remains unsolved, not that it does not exist. \citep[][p.~382]{zabell1992}
\end{quote}

\subsection{Is probability theory right for the job?}
\label{SS:false.conf}

That uncertainty quantification must be formulated using probability theory is almost exclusively taken for granted, at least in the statistics literature.  But it's worth asking if probabilistic uncertainty quantification is capable of achieving the reliability desiderata described in Section~\ref{SS:prob.uq}.  Spoiler alert---the answer to this question is {\em No}.  

My claim is that probabilistic uncertainty quantification is incapable of achieving the notion of reliability introduced in Section~\ref{SS:prob.uq}. For any mapping $z \mapsto \prior_z$ that represents a probabilistic quantification of uncertainty about the unknown $\Theta$, given $Z=z$, there are false hypotheses $H$, i.e., $H \not \ni \Theta$, such that $\prior_Z(H)$ tends to be large as a function of $Z \sim \prob_\Theta$.  If a large $\Pi_Z(H)$ is interpreted as ``confidence'' in the truthfulness of $H$, then the undesirable cases where $H$ is false and $\prior_Z(H)$ is large can be characterized as instances of {\em false confidence}, and it turns out that every form of probabilistic uncertainty quantification---Bayes, fiducial, etc.---is afflicted by it.  That is, however large a probability must be in order to instill confidence, there exists false hypotheses about which the posterior instills confidence with arbitrarily high rate, hence a risk of ``systematically misleading conclusions.''  The {\em false confidence theorem} below makes this precise.

% medical terminology that might work to use:
% Basu paper and BELIEF'24 paper help to diagnose false confidence
% can't cure false confidence in a non-trivial way, but it can be managed

\newcommand{\fcscore}{\text{\sc fcr}}

\begin{thm}[\citealt{balch.martin.ferson.2017}]
\label{thm:fct}
Let $z \mapsto \prior_z$ determine a data-dependent probability distribution absolutely continuous with respect to Lebesgue measure on $\TT$.  Define the false confidence rate associated with $z \mapsto \prior_z$ as 
\begin{equation}
\label{eq:fcscore}
\fcscore_\prior(\alpha, H) = \sup_{\theta \not\in H} \prob_\theta\{ \prior_Z(H) > 1-\alpha \}, \quad \alpha \in (0,1), \quad H \subset \TT. 
\end{equation}
Then for any $(\alpha, \tau) \in (0,1)^2$, there exists hypotheses $H$ such that $\fcscore_\prior(\alpha, H) > \tau$. 
\end{thm}

In other words, whatever level $1-\alpha$ is needed to instill confidence in the truthfulness of a hypothesis, there exists hypotheses $H$ for which the false confidence rate $\fcscore_\prior(\alpha, H)$ is arbitrarily large.  Since this is true for any probabilistic uncertainty quantification, the false confidence phenomenon can't be blamed on a poor choice of prior distribution, on issues underlying Fisher's fiducial argument, etc.  It's a shortcoming of probability theory with respect to data-driven uncertainty quantification about unknowns.  %Further discussion about the implications of the false confidence theorem will be discussed below.  

A criticism that must be dispelled right away is that high false confidence rate is only for ``trivial'' hypotheses, e.g., $H$ being the complement of a set with zero Lebesgue measure.  It's true that such ``trivial'' $H$ are afflicted by false confidence since $\prior_z(H) \equiv 1$ for all $z$, but false confidence occurs for lots of seemingly innocuous hypotheses, not just for these trivial ones. One of the many such examples is next.  
%Before moving on, it's important to mention that nothing in the problem setup that implies complement-of-a-``small''-set hypotheses should be ignored; after all, small hypotheses can be true and large hypotheses can be false, so size shouldn't influence the reliability of one's uncertainty quantification.  Fisher's familiar warnings about misinterpreting likelihood as probability are relevant here:
%\begin{quote} \label{quote:fisher.likelihood}
%%The function of the $\theta$'s... 
%The [likelihood function] is not however a probability and does not obey the laws of probability; it involves no differential element $d\theta_1 \, d\theta_2 \, d\theta_3$...; it does none the less afford a rational basis for preferring some values of $\theta$, or combination of values of the $\theta$'s, to others. \citep[][p.~552]{fisher1930}
%\end{quote} 
%He's clearly arguing that likelihood alone can be used to assess the compatibility of the observed data and a hypothesis $H$---a ``combination of values''---about $\Theta$.  If the likelihood has no differential element, then compatibility isn't assessed by integrating over $H$ and, therefore, the size of $H$ must not be relevant to these compatibility considerations.  

% Hypotheses that are obviously false aren't the ones practitioners are interested in testing.  Real discoveries will surely come down to edge cases, where the answer isn't obvious from looking at an appropriate plot, where the true $\Theta$ is near the boundary of $H$ and $H^c$, so it's the ``just barely false'' $H$'s that are of primary interest.  

\begin{figure}[t]
\begin{center}
\scalebox{0.5}{\includegraphics{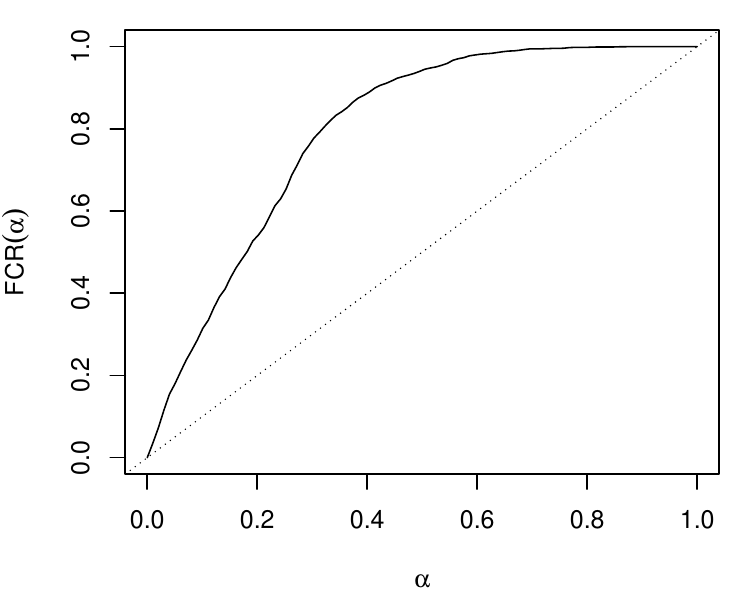}}
\end{center}
\caption{A lower bound on the (empirical) false confidence rate $\alpha \mapsto \fcscore_\prior(\alpha; H)$ of the Bayes posterior corresponding to the hypothesis $H$ stated in the text.}
\label{fig:fc}
\end{figure}

Consider a simple linear regression model, where $Z=(Z_1,\ldots,Z_n)$, with $Z_i = (X_i, Y_i)$ and $(Y_i \mid x_i, \theta) \ind \prob_{\theta, i} := \nm(\beta_0 + \beta_1 x_i, \sigma^2)$, for $\theta=(\beta_0, \beta_1, \sigma^2)$.  As is customary, I'll treat the $x_i$'s as given constants that are fixed; in the simulation below, I sampled these values independently from $\unif(-2,2)$. A standard Bayesian default prior corresponds to a certain choice of the conjugate normal--inverse gamma prior, so the corresponding posterior inference is rather straightforward.  To see how false confidence can creep in, suppose one is interested in the hypothesis 
\[ H = \{(\beta_0, \beta_1, \sigma^2): -\beta_0 / \beta_1 > -1\}, \]
which amounts to a hypothesis that the root of the regression function is larger than $-1$.  Suppose the true $\Theta$ is $(0.3, 0.1, 1)$, so that the above hypothesis is actually {\em false}.  Figure~\ref{fig:fc} shows a plot of (a lower bound on) the false confidence rate $\alpha \mapsto \fcscore(\alpha, H)$ associated with this Bayesian posterior distribution based on 1000 data sets of size $n=25$.  Note that even the lower bound is very high across the entire range of $\alpha$.  This tendency of the Bayesian posterior to assign relatively high probability to a false hypothesis is what creates a risk of systematically misleading conclusions.  

The risk associated with formulating statistical uncertainty quantification in terms of probability theory isn't new.  For example, \citet{fraser.cd.discuss} writes
\begin{quote}
[\citet{xie.singh.2012}]  
%Xie and Singh 
are thus recommending that we ignore the restriction to confidence sets or equivalent, and free confidence to allow the production of parameter distributions. Certainly distributions are easier to think about, are largely in accord with Fisher’s original proposal, and are more in the freedom of the Bayes approach, but they do overlook inherent risks... 
%as the preceding references indicate.} 
%\citep{fraser.cd.discuss}
\end{quote}
These risks mainly concern unreliability that results when marginalization is performed \citep[e.g.,][]{dawid.stone.zidek.1973, fraser2011, balch.martin.ferson.2017}.  When \citet{schweder.hjort.cd.discuss} warn that ``joint [confidence distributions] should not be sought, we think, since they might easily lead the statistician astray,'' they're concerned that users won't resist the temptation to do the familiar probabilistic marginalization, thereby creating a risk of unreliability.  What causes false confidence or this risk of unreliability more generally remains an open question.  The current conjecture is that false confidence tends to occur when the hypothesis concerns values of a non-linear function of the model parameters \citep{martin.belief2024}, like the $H$ above is a hypothesis about the ratio; the analyses in \citet{fraser2011} and \citet{fraser.etal.2016} offer similar take-away messages.  

\subsection{If not probability theory, then what?}
\label{SS:what}

The issues related to the unreliability of probabilistic uncertainty quantification in the context of statistical inference are general, not specific to any particular brand of probabilistic uncertainty quantification.  Therefore, to resolve these issues one must move beyond probabilistic uncertainty quantification to something else.  This ``something else'' should be probability-like, so that uncertainty quantification still makes sense, but it can't be additive.  Capacities, introduced in \citet{choquet1953}, are non-additive set functions, and I'll show below that a special type of data-dependent capacity can achieve the desired reliability properties that no ordinary probability can.  

Loosely speaking, {\em imprecise probabilities} are normalized capacities having additional properties making them suitable models for quantifying uncertainty.  Examples that might be familiar to statisticians include belief functions, first developed by \citet{dempster1966, dempster1967, dempster1968a} and later formalized by \citet{shafer1976}, which are infinitely-monotone capacities; the 2-monotone capacities that can be found in robustness investigations \citep[e.g.,][]{huber1973.capacity, huber1981, wasserman.kadane.1990, wasserman1990, berger1984}; and the so-called generalized Bayes framework of lower previsions advanced by \citet{walley1991}.  No specific details about these particular forms of imprecise probability will be needed in what follows.  Here I want to focus on {\em imprecision} itself and its role.  

In textbooks, ordinary or precise probability theory is presented in the context of a chance experiment---e.g., rolling an equally weighted six-sided die---where the specifics are all clearly stated but the outcome can't be predicted with certainty.  Then probability quantifies one's uncertainty about whether the unpredictable outcome of the experiment will satisfy such-and-such property.  This type of uncertainty is {\em aleatory}.  But what if the experiment's specifics weren't all clearly stated?  If there was ambiguity about the to-be-rolled die's configuration---e.g., {\em maybe} half of the faces are labeled with ``3,'' {\em maybe} it's asymmetrically weighted favoring ``6,'' etc.---then clearly there's no single probability that accurately captures uncertainty about the outcome.
%, since the different die configurations correspond to different degrees and forms of aleatory uncertainty, hence produce different numerical probability values.  
This ambiguity is an example of {\em epistemic} uncertainty, and this can't be accommodated using ordinary probability theory.  In the extreme case of ignorance about the die, applying the principle of indifference and assuming equal probabilities on each face is unacceptable: how can the assessments not differ between virtually orthogonal situations where one is ignorant about the die and where one is sure that the die is fair?  The issue isn't the assumption of fairness, it's the belief that a single probability can capture both aleatory and epistemic uncertainty.  Imprecise probability aims to directly address the epistemic uncertainty, the specification-ambiguity.  So, imprecision is not an inadequacy due to a poor assessment; it's about making an honest effort to faithfully capture all of what's uncertain.  

%which is related to the aphorism attributed to \citet{box1976}: ``All models are wrong.'' 

A bit more specifically, an imprecise probability corresponds to a lower and upper probability pair which, for our purposes here, is denoted by $(\lprob, \uprob)$.  That the two elements of the pair are related to each other will become clear shortly (see, also, Section~\ref{SS:construction}), but generally neither are probability measures.  The simplest interpretation of $(\lprob,\uprob)$ is literally as lower and upper bounds on a collection of precise probabilities.  For concreteness, while still taking some mathematical liberties for the sake of conceptual clarity, consider the die example above and let $\mathcal{D}$ denote the set of all the dice that could possibly be used in the experiment; this set captures the epistemic uncertainty.  Each die $D \in \mathcal{D}$ has an associated probability $\prob_D$ quantifying the aleatory uncertainty about the outcome of rolling die $D$.  To also capture the epistemic uncertainty, one can use 
\[ \lprob(\cdot) = \inf_{D \in \mathcal{D}} \prob_D(\cdot) \quad \text{and} \quad \uprob(\cdot) = \sup_{D \in \mathcal{D}} \prob_D(\cdot). \]
That $(\lprob,\uprob)$ are linked together should now be clear---they're both tied to the collection $\{\prob_D: D \in \mathcal{D}\}$.  In the extreme case where one is ignorant about the to-be-rolled die, the set $\mathcal{D}$ contains literally all possible dice and, consequently, $\lprob(A) = 0$ and $\uprob(A) = 1$ for all $A \not\in \{\varnothing, \varnothing^c\}$.  More generally, it's now clear that imprecise probabilities are more nuanced than ordinary or precise probabilities, and the motivation behind this added complexity is to properly handle both aleatory and epistemic uncertainty.  

This is relevant to the present goal of uncertainty quantification about $\Theta$ because, with {\em a priori} ignorance, epistemic uncertainty dominates.  From this perspective, the idea that data is sufficiently informative to justify a mapping from the vacuous imprecise prior---complete ignorance---to a posterior that's both fully precise and reliable is completely unrealistic.  {\em Imprecision is imperative}.  Indeed, the generalized Bayes rule of \citet{walley1991} applied to a fully vacuous prior returns a vacuous posterior, meaning it's impossible to learn (in a Bayesian way) when one is {\em a priori} ignorant; see, also, \citet{kyburg1987}, \citet{walley2002}, and, more recently, \citet{gong.meng.update}.  Non-Bayesian learning isn't susceptible to such criticism, but this upper hand apparently comes at the rather steep price of having to abandon probabilism altogether in favor of procedures with no natural fixed-data uncertainty quantification interpretation; see Zabell's quote on page~\pageref{quote:zabell}.  

My claim is, however, is that many of these non-Bayesian learning strategies do correspond to imprecise-probabilistic or, more specifically, {\em possibilistic} uncertainty quantification, it's just that no one has realized it.  Although imprecise probability theory didn't exist when Fisher was active, there are passages in his texts that suggests he may have anticipated an inexact or imprecise probability theory:
\begin{itemize}
%\item ``As Bayes perceived, the concept of mathematical probability affords a means, in some cases, of expressing inferences from observational data, involving a degree of uncertainty, and of expressing them rigorously...
%%in that the nature and extent of the uncertainty is specified with exactitude, 
%yet it is by no means axiomatic that the appropriate inferences...
%%, though in all cases involving uncertainty, 
%should always be rigorously expressible in terms of this same concept.'' (ibid., p.~40) % \citep[][p.~40]{fisher1973}
\item ``[A p-value] is more primitive, or elemental than, and does not justify, any exact probability statement about the proposition'' 
%(ibid., p.~46) 
\citep[][p.~46]{fisher1973}
\vspace{-2mm}
\item 
%``[Confidence regions] were I think developed and advocated under the impression that in a wider class of cases they could provide information similar to that of the probability statements derived by the fiducial argument.  
``It is clear, however, that no exact probability statements can be based on [confidence limits]'' 
%them.'' 
(ibid., p.~74) %\citep[][p.~74]{fisher1973}
\end{itemize} 
Presumably, non-Bayesians aren't opposed to fixed-data uncertainty quantification interpretations, they just don't know how to justify this without going a Bayesian route and possibly jeopardizing reliability.  The developments described below show how to get possibilistic uncertainty quantification that is both reliable and efficient.

\section{Possibilistic inferential models}
\label{S:details}

\subsection{Perspective}

An {\em inferential model} (IM) is a mapping $z \mapsto (\lPi_z, \uPi_z)$ from data $z$ to a lower and upper probability pair, $\lPi_z$ and $\uPi_z$, defined on $\TT$ to quantify uncertainty about the unknown $\Theta$, i.e., a model for uncertainty quantification and inference.  This includes the probabilistic approaches discussed earlier as a special case.  But \citet{imbasics, imbook} were focused on IMs that are provably reliable, which requires special care and consideration, so I'll follow their lead and not give too much attention to this general definition.  

Martin and Liu's original construction involved three steps.  The first step was an expression of the posited model via a data-generating process, or an {\em association}, say, $Z = a(\Theta,U)$, between observable data $Z$, unknown parameter $\Theta$, and an unobservable auxiliary variable $U$ with a known distribution.  The same association is used in Fisher's fiducial argument and in the generalizations proposed by \citet{dempster1967}, \citet{fraser1968}, and \citet{hannig.review}.  What's different about the IM approach is that attention focuses on the {\em unobserved value} $u$ of $U$, since the relation $z=a(\Theta,u)$ between the observed data $z$ and the unknown $\Theta$ must hold.  The second step, unique to Martin and Liu's IM construction, is where imprecision is introduced: they quantify uncertainty about the unobserved value $u$ using the distribution of a suitable random set.  Since the distribution of $U$ is known, it's easy to ensure that this uncertainty quantification about the unobserved $u$ is reliable.  Then the third step maps the random set for $u$ to a corresponding random set for $\Theta$, given $Z=z$, via the relation $z=a(\Theta,u)$, and the reliability offered by the former is immediately transferred to the latter.  Then uncertainty quantification about $\Theta$, given $Z=z$, is based on the distribution of that latter random set which, again, can be expressed as a belief function.  

This general approach is quite powerful, offering provably reliable, prior-free imprecise-probabilistic inference about $\Theta$, which was arguably Fisher's goal with the fiducial argument.  But the general IM approach also has some limitations.  First, it requires specification of an association ``$Z=a(\Theta,U)$,'' which is not uniquely determined by the model $\{\prob_\theta: \theta \in \TT\}$.  Second, in order for the inference to be efficient, non-trivial efforts to manipulate the association and reduce the dimension of the auxiliary variable are required \citep{imcond, immarg}.  Third, a justifiably ``optimal'' choice of the random set for the unobserved $u$ has remained elusive.  None of these on their own is a serious obstacle but, in aggregate, they create a barrier preventing the use of IMs, at least by non-experts, in practical applications.  A different, {\em generalized IM} approach was advanced in \citet{plausfn, gim} that specifically avoids the aforementioned limitations, and it's the recent developments around this version of the IM construction that I'll focus on here.  

One technical point concerning the random set-based construction is needed to provide background for what follows.  Theorem~4.3 in \citet{imbook} says the only admissible random sets for quantifying uncertainty about the unobserved $u$ are nested, i.e., for any two realizations of the random set, one is is a subset of the other.  While the distribution of random sets can generally be described by belief functions, the distribution of a nested random set corresponds to a special type of belief function, namely, a consonant belief function; see, e.g., \citet{shafer1976, shafer1987}.  Consonant belief functions correspond to {\em possibility measures} \citep[e.g.,][]{dubois.prade.book, dubois2006}, and these are similar to the probability distributions that statisticians are familiar with.  Since the aforementioned theorem implies that efficient IMs must take the form of possibility measures on $\TT$, I'll focus exclusively on {\em possibilistic IMs}; see, also, \citet{imposs}.  Readers unfamiliar with the basics of possibility theory should consult the brief Appendix~\ref{A:possibility} for the background relevant to the statistical developments below.

\subsection{Construction}
\label{SS:construction}

%{\color{red} Ranking and validification details, examples, pictures... two examples---exponential and gamma---that are neither trivial nor too familiar/simple... show how new approach is related to the old random set-based construction (maybe with details in the appendix)...}

The posited model $\{ \prob_\theta: \theta \in \TT \}$ and observed data $Z=z$ determine a likelihood $\theta \mapsto p_\theta(z)$ function and a corresponding relative likelihood
\[ 
R(z,\theta) = \frac{p_\theta(z)}{\sup_\vartheta p_\vartheta(z)}, \quad \theta \in \TT.
\]
I will assume throughout that the denominator is finite for almost all $z$.  An imprecise-probabilistic (in fact, possibilistic) interpretation can be given to the relative likelihood directly, and this has been extensively studied \citep[e.g.,][]{shafer1982, wasserman1990b, denoeux2006, denoeux2014}.  But the raw, relative likelihood-based possibilistic uncertainty quantification has issues similar to those discussed above for probabilistic uncertainty quantification---one apparently has no control on the false confidence rate.  But the relative likelihood serves an important role, namely, ranking different parameter values based on their compatibility with the observed data $Z=z$, which was the role that Fisher envisioned.  Arguably the relative likelihood $\theta \mapsto R(z,\theta)$ above is the ``best'' such ranking function since it's a minimial sufficient statistic; see, also, Remark~\ref{re:rellik} in Appendix~\ref{A:remarks}.  But this is not the only ranking function one might entertain; see Sections~\ref{S:marginal}--\ref{S:beyond}.  

\begin{figure}[t]
\begin{center}
\scalebox{0.5}{\includegraphics{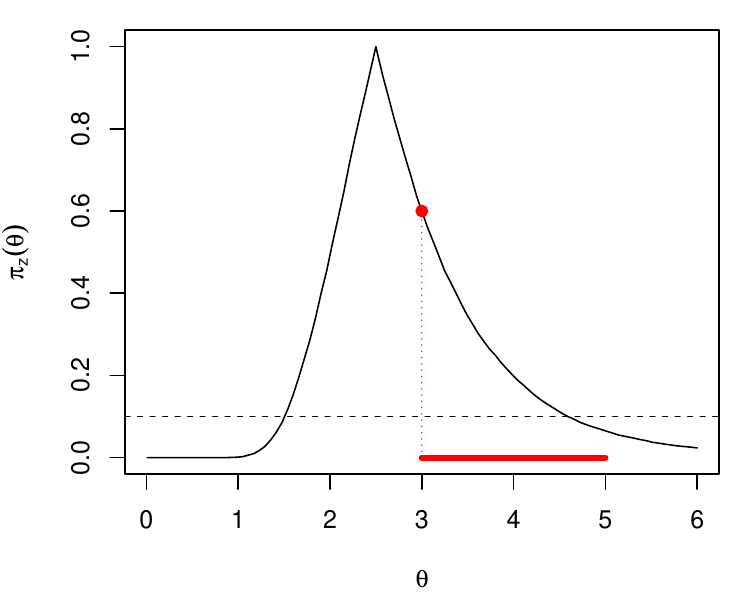}}
\end{center}
\caption{Illustration of a possibility contour $\theta \mapsto \pi_z(\theta)$ and how the possibility measure $\uPi_z$ is determined from it.  Here, the hypothesis $H$ is the interval $[3,5]$ and the maximum value $\uPi_z(H)$ on contour is highlighted. Horizontal dashed line at $\alpha=0.1$ determines the 90\% confidence interval $C_{0.1}(z)$ in \eqref{eq:region}.}
\label{fig:pl.toy}
\end{figure}

The second step of the possibilistic IM construction is ``validifying'' \citep{martin.partial} the relative likelihood (or other ranking function).  This amounts to applying a version of the {\em probability-to-possibility transform} \citep[e.g.,][]{dubois.etal.2004, hose2022thesis}, and it returns the possibilistic IM's contour function:
\begin{equation}
\label{eq:contour}
\pi_z(\theta) = \prob_\theta\bigl\{ R(Z,\theta) \leq R(z, \theta) \bigr\}, \quad \theta \in \TT.
\end{equation}
This is indeed a possibility contour since the function is (clearly) non-negative and satisfies $\sup_\theta \pi_z(\theta) = 1$; in fact, the supremum is attained at the maximum likelihood estimator, $\hat\theta_z = \arg\max_\theta R(z,\theta)$, i.e., $\pi_z(\hat\theta_z)=1$.  The corresponding possibility measure, or upper probability, is defined via optimization as 
\begin{equation}
\label{eq:maxitive}
\uPi_{z}(H) = \sup_{\theta \in H} \pi_{z}(\theta), \quad H \subseteq \TT.
\end{equation}
%Note that this optimization-based formulation involves no differential element and, therefore, the size of $H$ has no bearing at all on $\uPi_z(H)$; see the Fisher quote on page \pageref{quote:fisher.likelihood}.  
The corresponding necessity measure, or lower probability, is defined via conjugacy: $\lPi_{z}(H) = 1 - \uPi_{z}(H^c)$.  An example of $\pi_z$ and its corresponding $\uPi_z$ is shown in Figure~\ref{fig:pl.toy}.  An explanation of how this new possibilistic IM construction is related to the original random set-based IM construction is given in Appendix~\ref{A:old.vs.new}. 

More details follow, but the key point is that $H \mapsto \{\lPi_z(H), \uPi_z(H)\}$ reliably quantifies uncertainty about $\Theta$, given $Z=z$.  For example, $H$ such that $\uPi_z(H)$ is small are implausible because the data provides little/no support in favor of the truthfulness of ``$\Theta \in H$'' and, similarly, $H$ for which $\lPi_z(H)$ is large are not only plausible but believable since the data provides strong support for the truthfulness of ``$\Theta \in H$.''

\subsection{Sampling-based reliability properties}
\label{SS:properties}

\subsubsection{Validity}
\label{SSS:valid}

% Walley (2002?)~strikes a middleground, but this too has a certain form of dilation---the inference is inefficient compared to the $n^{-1/2}$-consistent solutions in textbooks, etc. 

The essential reliability property of the possibilistic IM is {\em strong validity}.

\begin{thm}
\label{thm:valid}
The possibilistic IM is strongly valid in the sense that 
\begin{equation}
\label{eq:valid}
\sup_{\theta \in \TT} \prob_\theta\bigl\{ \pi_{Z}(\theta) \leq \alpha \bigr\} \leq \alpha, \quad \text{$\alpha \in [0,1]$}. 
\end{equation}
\end{thm}

This corresponds to the familiar result for p-values and is a direct consequence of the probability integral transformation taught in the basic math-stat course sequence.  While this and some---but not all---of the results below might be familiar in the context of p-values, it's important to keep in mind that p-values are commonly used for isolated significance tests, not as building blocks for a broad, mathematically rigorous framework for reliable uncertainty quantification.  Moreover, the results here align with familiar p-value considerations only because I'm assuming vacuous prior information; the more general case is covered in \citet{martin.partial2} and discussed briefly in Section~\ref{S:discuss} below.  

Strong validity has a number of important consequences.  First, \eqref{eq:valid} immediately implies that the upper $\alpha$ level set of the possibility contour is a $100(1-\alpha)$\% confidence region.  Note that Bayesian and fiducial credible sets can generally only reach confidence set status asymptotically, as sample size approaches infinity.   Figure~\ref{fig:pl.toy} gives a depiction of the upper level set $C_\alpha(z)$, with $\alpha=0.1$.  

\begin{cor}
\label{cor:coverage}
The upper $\alpha$ level set of the possibilistic IM contour,
\begin{equation}
\label{eq:region}
C_\alpha(z) = \{\theta \in \TT: \pi_{z}(\theta) \geq \alpha\}, \quad \alpha \in [0,1],
\end{equation}
is a $100(1-\alpha)$\% confidence set, i.e., $\sup_{\theta \in \TT} \prob_\theta\bigl\{ C_\alpha(Z) \not\ni \theta \bigr\} \leq \alpha$.  
\end{cor}

Second, a more basic property called {\em validity}, introduced in \citet{imbasics}, is naturally a consequence of strong validity and the relation \eqref{eq:maxitive} between the possibility contour and the possibility measure.  The possibilistic IM is said to be {\em valid} if 
\begin{equation}
\label{eq:valid.alt}
\sup_{\theta \in H} \prob_\theta\bigl\{ \uPi_Z(H) \leq \alpha \bigr\} \leq \alpha, \quad \text{for all $\alpha \in [0,1]$, $H \subseteq \TT$}.
\end{equation}
In words, property \eqref{eq:valid.alt} implies that it's a small probability event that the IM assigns small upper probability to a true hypothesis.  Thanks to the relationship between the IM's lower and upper probability, this can equivalently be expressed as 
\[ \sup_{\theta \not\in H} \prob_\theta\bigl\{ \lPi_Z(H) > 1-\alpha \bigr\} \leq \alpha, \quad \text{for all $\alpha \in [0,1]$, $H \subseteq \TT$}. \]
Then, just like the interpretation of \eqref{eq:valid.alt}: it's a small probability event that the IM assigns large lower probability to a false hypothesis.  

A natural question is why should the same quantity, $\alpha$, appear both inside and outside of the curly brackets in the above two displays.  The reason is that the interpretation of numerical probabilities is context-independent.  That is, while what probabilities are deemed ``small'' and ``large'' may vary from one individual to the next, a statement like {\em the probability is 0.1} means the same thing to a particular individual whether the topic is tomorrow's weather or the reliability of my data-driven uncertainty quantification about an unknown $\Theta$.  So, the scale on which the probabilities about $\Theta$ would be interpreted is exactly the same as that on which the model-based probabilities about $Z$ would be interpreted.  Hence, the same $\alpha$---representing whatever value one interprets as ``small''---appears inside and outside of the probability statement in \eqref{eq:valid.alt}.

Validity in the sense of \eqref{eq:valid.alt} has several important consequences concerning the construction of hypothesis testing procedures and the control of false confidence rate.  Concerning the latter, while there's no way to eliminate false confidence altogether, except by refusing to learn, i.e., by setting $\lPi_z \equiv 0$ for all $z$, it's important to control $\fcscore$ in terms of the confidence threshold $\alpha$.  Bayesian and fiducial uncertainty quantification offers no control over the $\fcscore$, but possibilistic IMs do. 

\begin{cor}
The possibilistic IM is valid in the sense of \eqref{eq:valid.alt}.  Consequently: 
\begin{itemize}
\item The test ``reject $H$ if $\uPi_z(H) \leq \alpha$'' controls the Type~I error probability at level $\alpha$.  
\vspace{-2mm}
\item If $\fcscore_{\lPi}(\alpha, H)$ is the false confidence rate associated with the possibilistic IM's lower probability $z \mapsto \lPi_z$, analogous to \eqref{eq:fcscore}, then $\fcscore_{\lPi}(\alpha, H) \leq \alpha$ for all $H \subset \TT$. 
\end{itemize} 
\end{cor}

Finally, while validity and strong validity were more-or-less treated as equivalent properties in, e.g., \citet{imbasics, imbook}, it's important to emphasize that strong validity \eqref{eq:valid} is indeed stronger than validity \eqref{eq:valid.alt}.  This was first established in \citet{cella.martin.probing}, where strong validity in \eqref{eq:valid} and a uniform-in-hypotheses version of \eqref{eq:valid.alt} were shown to be equivalent.  See Remark~\ref{re:uniform} in Appendix~\ref{A:remarks} for further explanation. 

%I'll skip these details for the sake of space.

\subsubsection{Efficiency}
\label{SSS:efficiency}

Validity is not the only reliability-relevant concept.  In fact, validity is incredibly easy to achieve---just take $\lPi_z(H) \equiv 0$ and $\uPi_z(H) \equiv 1$ for all $z$ and all (non-trivial) $H$.  The problem, of course, with this latter choice is that the output exactly agrees with the vacuous prior; nothing is learned from data.  So the goal is to get the most out of the available data but in such a way that reliability is preserved.  I have already referred and will continue to refer to this complementary notion of reliability as {\em efficiency}.  The issues here are fundamental, related to basic conceptual notions developed in the early 19th century by Legendre and Gauss concerning least squares and, more generally, the combination of observations \citep[e.g.,][]{stigler1986} and later, in the 20th century, sufficient statistics, Fisher information, the Cram\'er--Rao lower bound, etc.  

Carefully combining information across distinct sources is critical for efficient inference.  In the early IM developments, \citet{imcond} handled this combination manually, by suitably manipulating the association that links data $Z$, parameter $\Theta$, and auxiliary variable $U$.  They ``rediscovered'' classical dimension reduction techniques, such as sufficiency and conditioning on ancillary statistics; they also developed some new insights, beyond the scope of this review.  While their manual approach offers greater flexibility and, in turn, the potential for greater efficiency, this is often difficult to carry out.  The present construction, with relative likelihood-based ranking, automatically combines observations in a way that's ``optimal,'' at least in certain cases, without any manual effort.  The results reviewed below, from \citet{imbvm.ext}, show that the valid possibilistic IM above is asymptotically efficient in a familiar sense.  Consequently, the IM's exact validity (via imprecision) comes at no cost in terms of efficiency.  

What's summarized below is a possibility-theoretic version of the celebrated {\em Bernstein--von Mises theorem} that appears in the Bayesian and (generalized) fiducial literature, ensuring that the output is asymptotically Gaussian with covariance matrix that agrees with the Cram\'er--Rao lower bound.  The Bernstein--von Mises theorem is important for Bayes and fiducial because it ensures that the credible sets are asymptotic confidence sets.  For the possibilistic IM, the contour's level sets are automatically confidence sets (Corollary~\ref{cor:coverage}), so Theorem~\ref{thm:bvm} below strictly concerns the IM's asymptotic efficiency.  

To make sense out of this, I need to introduce the notion of a Gaussian possibility measure.  Let $g_v$ denote the $d$-dimensional Gaussian probability density function, parametrized by a covariance matrix $v \in \RR_+^{d \times d}$.  Define the corresponding Gaussian possibility measure as the outer possibilistic approximation \citep[e.g.,][and Appendix~\ref{A:possibility}]{dubois.prade.1990} of $\nm_d(0, v)$.  
%, i.e., the possibility measure with the smallest credal set that contains $\nm_D(m, v)$.  
%Analogous to the formula \eqref{eq:contour} presented in Section~\ref{S:background}
Equivalently, the Gaussian possibility measure is just the probability-to-possibility transform of a Gaussian distribution, with contour 
\begin{align*}
\gamma_{v}(y) = \prob\{ g_v(Y) \leq g_v(y) \} 
%& = \prob\{ e^{-\frac12 (Y-\mu)^\top \sigma^{-1} (Y - \mu)} \leq e^{-\frac12 (y-\mu)^\top \sigma^{-1} (y - \mu)} \} \\
%& = \prob\{ (Y-\mu)^\top \sigma^{-1} (Y - \mu) \geq (y-\mu)^\top \sigma^{-1} (y - \mu) \} \\
= 1 - F_d( y^\top v^{-1} y ), \quad y \in \RR^d, \quad Y \sim \nm_d(0, v), 
\end{align*}
where $F_d$ is the $\chisq(d)$ distribution function.  
%Write $\uGamma_{m,v}$ for the corresponding possibility measure obtained by maximizing the contour as described above.  
Since we're concerned with large-sample approximations, replace the generic $Z$ by a vector $Z^n = (Z_1,\ldots,Z_n)$ of iid components, and write $\pi_{Z^n}$ for the IM's possibility contour.  The general theorem in \citet{imbvm.ext} says that, under conditions comparable to those in  \citet{lecam1970}, $\pi_{Z^n}$ and a suitable Gaussian possibility contour---one whose variance agrees with the Cram\'er--Rao lower bound---merge with probability converging to 1 as $n \to \infty$.  Below is a simpler-to-state version that assumes the stronger conditions typically found in textbooks. 

\begin{thm}[\citealt{imbvm.ext}---simpler version]
\label{thm:bvm}
Under the usual Cram\'er conditions, 
%presented in Appendix~\ref{A:large.sample}, 
the possibilistic IM's contour $\pi_{Z^n}$ satisfies 
%\[ \sup_{\theta \in \mathcal{T}} \bigl| \pi_{Z^n}(\theta) - \gamma_{Z^n}(\theta) \bigr| \to 0 \quad \text{in $\prob_\Theta$-probability as $n \to \infty$}, \]
%where $\mathcal{T}$ is an arbitrary compact subset of $\TT \subseteq \RR^d$, $\gamma_{Z^n}(\theta) = \gamma_{\Theta + n^{-1/2}\Delta_\Theta(Z^n), (n I_\Theta)^{-1}}(\theta)$ is the centered and scaled Gaussian possibility contour, $I_\Theta$ the Fisher information, and $\Delta_\Theta(X^n)$ a sequence of random variables that converges in distribution to $\nm_d(0, I_\Theta^{-1})$.  If, in addition, the likelihood function is twice differentiable, then 
\[ \sup_{u \in \mathcal{U}} \bigl| \pi_{Z^n}(\hat\theta_{Z^n} + J_{Z^n}^{-1/2} u) - \gamma(u) \bigr| \to 0 \quad \text{in $\prob_\Theta$-probability as $n \to \infty$}, \]
where $\mathcal{U}$ is an arbitrary compact subset of $\RR^d$, $\hat\theta_{Z^n}$ is the maximum likelihood estimator, $J_{Z^n}$ is the observed Fisher information matrix, and $\gamma$ is the standard Gaussian possibility with identity covariance matrix.
\end{thm}

%The reader of course will recognize that the covariance matrix in the Gaussian possibility measure agrees with that which appears in the Cram\'er--Rao lower bound.  It's in this sense that the IM output is asymptotically efficient.  

For further details, discussion, and illustrations, see \citet{imbvm.ext}.  To summarize, at least when prior information about the model parameter $\Theta$ is vacuous and there's no explicit priority for one assertion about $\Theta$ versus another, the relative likelihood-based possibilistic IM is both exactly valid and asymptotically optimal in the classical sense of matching the Cram\'er--Rao lower bound.  But beyond this standard case, improvements may be possible, and that's the focus of Sections~\ref{S:marginal}--\ref{S:beyond} below.

\subsection{Conditional, fixed-data properties}
\label{SS:conditional}

While it's natural and important to consider the sampling-based properties of the IM output, an often overlooked advantage of (imprecise) probabilistic uncertainty quantification is that it offers a fully conditional, fixed-data interpretation.  This angle isn't often treated in the literature on default-prior Bayes, (generalized) fiducial, IMs, etc.  Unfortunately, I also don't have the space here for a careful treatment, but see Appendix~\ref{A:conditional}.

\subsection{Computation}
\label{SS:computation}

Until recently, only naive and relatively inefficient strategies for computing the IM contour were available.  In particular, the go-to strategy was to approximate $\pi_z$ via 
\begin{equation}
\label{eq:pi.naive}
\pi_{z}(\theta) \approx \frac1M \sum_{m=1}^M 1\{ R(Z_{m,\theta}, \theta) \leq R(z, \theta) \}, \quad \theta \in \TT,
\end{equation}
where $Z_{m,\theta}$ are independent copies of data $Z$, drawn from $\prob_\theta$, for $m=1,\ldots,M$.  The above computation is feasible at a few different $\theta$ values, but all the practically relevant calculations---e.g., identifying the confidence set in \eqref{eq:region}---require evaluation over a fine grid covering the relevant portions of $\TT$, and this is expensive.  Some improvements to the basic strategy above were presented in \citet{hose.hanss.martin.belief2022}.  Strategies based on asymptotic approximations are also available, and discussed briefly in Remark~\ref{re:comp} in Appendix~\ref{A:remarks}.  

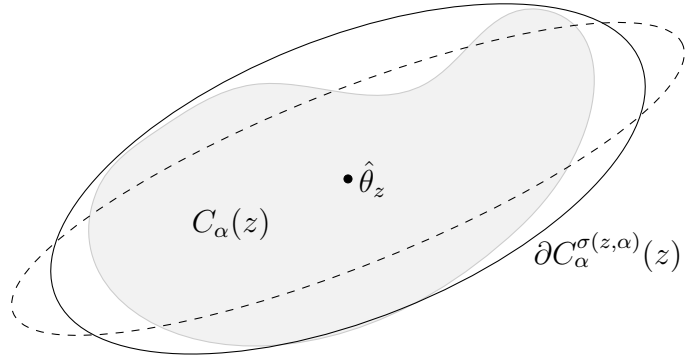
\begin{figure}[t]
\begin{center}
\scalebox{0.8}{
\begin{tikzpicture}
\filldraw[color=black!20, fill=black!5] (-1.25, 1.25) to[closed, curve through = { (1, 1.25) (2.5, 2.25) (2.25, -0.65) (-3.3, -0.15) (-2.5, 0.65) }] (-1.25, 1.25); 
\filldraw[color=black, fill=black] (0,0) circle (1.5pt) node[anchor=west]{$\hat\theta_z$};
\draw (-2.2, -0.6) circle (0pt) node[anchor=west]{$C_\alpha(z)$}; 
%\filldraw[dashed, rotate=50] (0.5,0) ellipse (2.1 and 1.1);
\draw[rotate=22] (0,0) ellipse (4.15 and 1.86);
\draw[dashed, rotate=22] (0,0) ellipse (4.75 and 1.15);
%\draw[dashed, rotate=22] (0,0) ellipse (2.75 and 2.5);
%\filldraw[color=black, fill=black] (-1.25,1.25) circle (1.5pt);
%\filldraw[color=black, fill=black] (1, 1.25) circle (1.5pt);
%\filldraw[color=black, fill=black] (2.5, 2.25) circle (1.5pt);
%\filldraw[color=black, fill=black] (-3, -1.15) circle (1.5pt);
\draw (2.30, -0.97) circle (0pt) node[anchor=west]{$\partial C_\alpha^{\sigma(z,\alpha)}(z)$}; 
%\filldraw[color=black, fill=black] (-2.5, 0.5) circle (1.5pt);
\end{tikzpicture}
}
\end{center}
\caption{Approximating the IM contour's $\alpha$-level set $C_\alpha(z)$ by an ellipse $C_\alpha^\sigma(z)$ with a ``good'' choice of $\sigma$ (solid) and with a ``bad'' choice of $\sigma$ (dashed).}
\label{fig:bean}
\end{figure}

\citet{immc} recently developed a new and efficient IM computational strategy.  That proposal replaces (most of) the naive contour evaluations in \eqref{eq:pi.naive} with Monte Carlo sampling from a sort of ``posterior distribution'' derived specifically from the IM output---not via Bayes's theorem.  The starting point for these developments is what's call the {\em credal set} \citep[e.g.,][Ch.~5]{levi1980} associated with the IM output.  In general, the credal set is just the collection of precise probabilities that are dominated by a given upper probability; in our present notation, this is defined as 
\begin{equation}
\label{eq:general.credal}
\cred(\uPi_z) = \{ \prior_z \in \text{probs}(\TT): \prior_z(H) \leq \uPi_z(H) \text{ for all measurable $H$} \}, 
\end{equation}
where $\text{probs}(\TT)$ is the set of all probability measures supported on (the Borel $\sigma$-algebra of) the parameter space $\TT$.  Since $\uPi_z$ is a possibility measure, there are simpler characterizations of the elements of $\cred(\uPi_z)$ in \eqref{eq:general.credal}.  Indeed, a now-classical characterization \citep[e.g.,][]{cuoso.etal.2001} of a possibilistic credal set is
\begin{equation}
\label{eq:poss.credal}
\prior_z \in \cred(\uPi_z) \iff \prior_z\{ C_\alpha(z) \} \geq 1-\alpha, \quad \text{all $\alpha \in [0,1]$}, 
\end{equation}
where $C_\alpha(z) = \{\theta: \pi_z(\theta) > \alpha\}$ is defined in \eqref{eq:region}.  Since $\{C_\alpha(z): \alpha \in [0,1]\}$ is a collection of confidence sets and since each element $\prior_z$ in $\cred(\uPi_z)$ assigns probability at least $1-\alpha$ to the $100(1-\alpha)$\% confidence set, it's fair to describe the elements of $\cred(\uPi_z)$ as {\em confidence distributions}.  Among these confidence distributions, clearly the ``best'' approximation $\prior_z^\star$ of $\uPi_z$ is one that achieves equality in the far-right inequality in \eqref{eq:poss.credal}, i.e., 
\begin{equation}
\label{eq:imcd}
\prior_z^\star\{ C_\alpha(z) \} = 1-\alpha, \quad \alpha \in [0,1]. 
\end{equation} 
From the information provided so far, it's not at all clear what form the elements of $\cred(\uPi_z)$ have and how the equality in \eqref{eq:imcd} might be achieved.  This leads to a more specific characterization, given in \citet{immc}, which is similar to but different from those in  \citet{wasserman1990} and in \citet{hose2022thesis}: $\prior_z \in \cred(\uPi_z)$ if and only if 
\begin{equation}
\label{eq:char}
\prior_z(\cdot) = \int_0^1 \kernel_z^\beta(\cdot) \, \marg_z(d\beta), 
\end{equation}
for some collection of probability distributions $\{\kernel_z^\beta: \beta \in [0,1]\}$ such that $\kernel_z^\beta$ fully supported on $C_\beta(z)$ for each $\beta \in [0,1]$, and some probability measure $\marg_z$ on $[0,1]$ such that a random variable with distribution $\marg_z$ is stochastically no smaller than $\unif(0,1)$.  From this characterization, it's easy to check that a distribution $\prior_z^\star \in \cred(\uPi_z)$ satisfies \eqref{eq:imcd} if and only if $\marg_z$ is $\unif(0,1)$ and $\kernel_z^\beta$ is fully supported on the boundary, $\partial C_\beta(z)$, of the confidence sets $C_\beta(z)$ for each $\beta \in [0,1]$.  If such a $\prior_z^\star$ exists, then it's called the {\em inner probabilistic approximation} of $\uPi_z$; see, also, \citet{dubois.prade.1990}.  Existence is typically not an issue \citep[see][for these details]{immc}, but the inner probabilistic approximation is not unique.  In the spirit of Keynes's principle of indifference, \citet{immc} recommends taking $\kernel_z^\beta$ to be $\unif\{ \partial C_\beta(z) \}$.  To summarize, a draw $(\Theta \mid z) \sim \prior_z^\star$, the inner probabilistic approximation, can be obtained by follow these two steps:
\begin{enumerate}
\item sample $A \sim \unif(0,1)$;
\vspace{-2mm}
\item sample $(\Theta \mid z, A=\alpha) \sim \unif\{ \partial C_\alpha(z) \}$.
\end{enumerate}
See the gamma example below and those in Appendix~\ref{A:examples}.  The second step is simple to state but non-trivial to carry out, in part because the boundary $\partial C_\alpha(z)$ isn't readily available.  \citet{immc} proposed a simple work-around, which is to bound the level set $C_\alpha(z)$ by a suitable ellipsoid, say, $C_\alpha^\sigma(z)$, indexed by a vector $\sigma$ of scaling parameters, in the sense that $C_\alpha(z) \subseteq C_\alpha^\sigma(z)$ for some $\sigma$.  In light of the asymptotic normality established in Theorem~\ref{thm:bvm}, the level set $C_\alpha(z)$ should be approximately ellipsoidal, so there should be no trouble finding an ellipsoid $C_\alpha^\sigma(z)$ in which $C_\alpha(z)$ fits snugly.  The appropriate choice of $\sigma$ would generally depend on both $z$ and $\alpha$, so I'll write it as $\sigma(z,\alpha)$.  See Figure~\ref{fig:bean} for an illustration.  The key point is that sampling from the boundary of an ellipsoid is easy, so if the ellipsoid is a good, not-anti-conservative approximation, then Step~2 can be reasonably approximated without sacrificing validity.  Details about the choice of $\sigma$ can be found in \citet{imvar.ext} and \citet{immc}.  Finally, the samples of $\Theta$ drawn according to the above procedure---which are ``probabilistic'' in the sense that they can be used to approximate probabilities---are  transformed back into a possibility measure using an empirical version of the probability-to-possibility transform described in Appendix~\ref{A:possibility}; for full details, see \citet[][Sec.~3.4]{immc}.

%The details of this new computational strategy are too intricate to present here, but the aforementioned IM-based ``posterior distribution,'' its origins, and its properties will be described in Section~\ref{S:reimagined}.  For the sake of completeness, details of this new computational strategy---which leans on the approximate normality coming from Theorem~\ref{thm:bvm}---are presented in Appendix~\ref{A:computation}. 

\subsection{Examples}
\label{SS:examples1}

\begin{figure}[t]
\begin{center}
%\scalebox{0.7}{\includegraphics{hamada_joint_pl}}
\subfigure[Normal example]{\scalebox{0.5}{\includegraphics{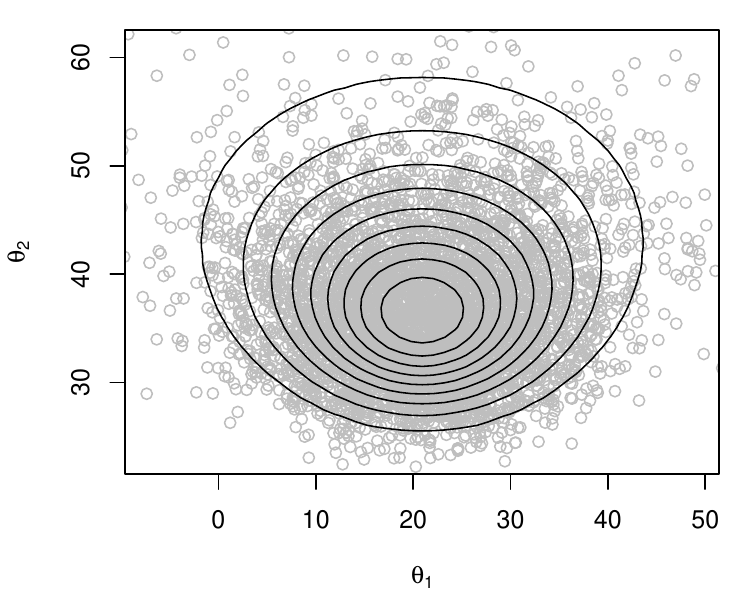}}}
\subfigure[Gamma example]{\scalebox{0.5}{\includegraphics{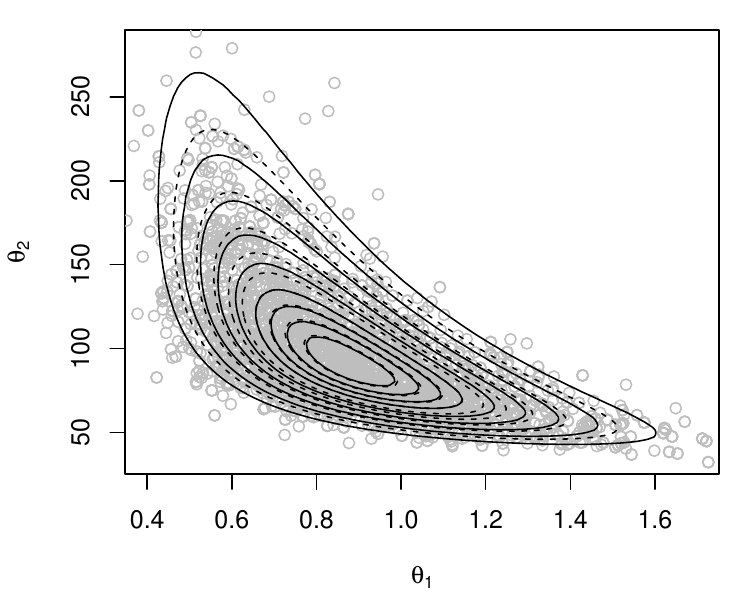}}}
\end{center}
\caption{Plot of the joint IM contour for $\Theta=(\Theta_1,\Theta_2)$ in the two examples presented in Section~\ref{SS:examples1}.  In both plots, the solid line is the possibilistic IM contour and the gray points are samples from the Jeffreys prior Bayes posterior. Dashed line in Panel~(b) is explained in the text.}
\label{fig:example1}
\end{figure}

\begin{gaussex}
The particular example considered here--famous for several reasons---is an analysis of Charles Darwin's data on the difference between cross- and self-fertilized plants.  This is the same example that appears in \citet[][Ch.~3]{fisher.doe}.  Context of the experiment aside, the data $z$ consists of $n=15$ differences in plant heights; the units are eighths of inches.  For simplicity, I'll model these as iid $\prob_\Theta = \nm(\Theta_1, \Theta_2^2)$, i.e., as normal random variables with unknown mean and standard deviation, $\Theta_1$ and $\Theta_2$, respectively.  The maximum likelihood estimators are $\hat\theta_{z,1} = 20.93$ and $\hat\theta_{z,2}=36.46$, respectively, and the relative likelihood is 
\[ R(z,\theta) = \Bigl( \frac{\hat\theta_{z,2}^2}{\theta_2^2} \Bigr)^{n/2} \exp\Bigl\{ -\frac{n(\hat\theta_{z,1} - \theta_1)^2}{2\theta_2^2} - \frac{n}{2} \Bigl( \frac{\hat\theta_{z,2}^2}{\theta_2^2} - 1 \Bigr) \Bigr\}, \quad \theta_1 \in \RR, \quad \theta_2 > 0. \]
Note that $R(Z,\theta)$ is a pivot when $Z$ consists of iid samples from $\prob_\theta$.  Although the distribution doesn't have a name or a simple form, it is easy to simulate and, therefore, the validification step can easily be done via Monte Carlo.  A plot of the possibility contour $\pi_z$ for Darwin's data is shown in Figure~\ref{fig:example1}(a).  For comparison, samples from the Bayesian posterior distribution based on Jeffreys prior is shown in the background in gray.  Both are pointing the same direction, but only the possibilistic IM offers reliable uncertainty quantification guarantees.  
\end{gaussex}

\begin{gammaex}
Consider the systems reliability application in \citet{hamada.etal.2004} that involves the breakdown times of $n=20$ machines.  At the 5\% significance level, the Kolmogorov--Smirnov test cannot reject the null hypothesis that these data are gamma, so I'll construct a possibilistic IM for $\Theta=(\Theta_1, \Theta_2)$, the shape and scale parameters of a gamma model.  It is easy to sample from a gamma distribution, and also easy to evaluate the relative likelihood, so the naive approximation \eqref{eq:pi.naive} of the IM contour is likewise easy; but carrying out these ``easy'' calculations over a sufficiently fine grid in the $(\theta_1,\theta_2)$-space is rather expensive, but the strategy in \citet{immc} discussed in Section~\ref{SS:computation} is more efficient.  A contour plot of this approximation of $\pi_z$ is shown in Figure~\ref{fig:example1}(b).  For comparison, I also show the large-sample Gaussian approximation (dashed line) and samples from a Bayesian posterior distribution using Jeffreys prior (gray points).  All three solutions point to roughly the same regions in the $(\theta_1,\theta_2)$-space that plausibly contain $\Theta$, and the shapes are similar.  Again, only the solid curves are backed by theory that guarantees reliability of uncertainty quantification. 
\end{gammaex}

\section{Implications for frequentists and Bayesians}
\label{S:reimagined}

%Zabell (page~\pageref{quote:zabell}) describes the fiducial argument as an attempt navigate safely between the dangerous Scylla and Charybdis corresponding to the unconditional, non-probabilistic frequentist perspective and the fully conditional, probabilistic Bayesian perspective, respectively.  Most authors writing about fiducial-related topics opt to play it safe, describing their proposal as a {\em unification} of the frequentist and Bayesian extremes.  But this unification argument doesn't challenge the status quo; in fact, unification gives everyone permission to ignore disagreements and inconsistencies and continue as they are.  While the status quo might be comfortable, it's apparently unhealthy for the discipline since we're still worried about ``missing the data science boat'' (Section~\ref{S:intro}).  Therefore, I think a less-safe approach might be worth taking here.  

\subsection{For frequentists} 
\label{SS:frequentists}

There are good reasons for frequentists to abandon probabilism.  The list of reasons includes the reliability-related warnings given in Section~\ref{SS:false.conf} above, foundational issues as detailed in \citet{mayo.book.2018}, and practically relevant points concerning probabilism's lack of flexibility, e.g., ``The idea that statistical problems do not have to be solved as one coherent whole is anathema to Bayesians but is liberating for frequentists'' \citep{wasserman.quote}.  But abandoning formal uncertainty quantification altogether because of the shortcomings associated with probabilism is extreme---like throwing the baby out with the bathwater.  Indeed, this abandonment is both unnecessary and harmful.  

I'll address the ``unnecessary'' claim first.  Frequentists have their favorite solutions to classical problems, so naturally they're reluctant to consider new frameworks that propose different solutions to these classical problems.  But the likelihood-based possibilistic IM described above often exactly agrees with the classical solution (modulo proper marginalization, if necessary; see Section~\ref{S:marginal}).  Moreover, the likelihood-based formulation can easily be generalized (Section~\ref{S:beyond}), providing more flexibility and broader agreement between the possibilistic IM solutions and go-to frequentist solutions.  Even more generally, the result formally stated and proved in Appendix~\ref{AA:char} says roughly the following: {\em for any test or confidence procedure concerning any feature $\Phi=f(\Theta)$ of the full parameter $\Theta$ that has frequentist error rate guarantees, there exists a valid possibilistic IM for $\Theta$---full uncertainty quantification!---that yields a test/confidence procedure for $\Phi$ that is at least as good as the given procedure}.  This result, which generalizes similar results in \citet{impval} and \citet{imchar}, has an important consequence: there are no genuine frequentist solutions---including those in classical textbooks and those that haven't been conceived of yet---that are out of the possibilistic IM framework's reach.  Therefore, frequentists are already using possibilistic IMs, so nothing about the proposed brand of uncertainty quantification can be objectionable to them; but they aren't taking full advantage of all that possibilistic IMs have to offer, which I'll discuss next.  

%\begin{itemize}
%\item Section~\ref{SS:properties} already proves that frequentist reliability priorities can easily be met within a formal, possibility theory-based uncertainty quantification framework. 
%\item In a nutshell, a thesis of \citet{mayo.book.2018} is that (the ``NHST'' caricature of) frequentism doesn't meet the needs of scientists seeking to reliably probe for hypotheses that are supported by data.  Data surely has more to offer than a binary accept/reject decision, but without a formal model for uncertainty quantification, she can only supplement the frequentist Type~I error rate and p-value with another measure that facilitates probing.  But the possibilistic IM output already features two components, a lower and upper probability, and \citet{cella.martin.probing} explain how the IM's lower probability plays the role of this support-probing measure.  The uniform validity proved in Corollary~\ref{cor:uniform} ensures that the IM's probing is reliable.  
%\end{itemize} 
That frequentists abandoning formal uncertainty quantification is harmful has been widely documented; several recent issues of {\em The American Statistician} have been devoted to this.
%, including the one prefaced by an editorial \citep{wasserstein.lazar.asa} that was misinterpreted as an official statement on p-values and statistical significance by the American Statistical Association, thereby prompting nomination of special task force to prepare and make an official statement.\footnote{\url{https://magazine.amstat.org/blog/2021/08/01/task-force-statement-p-value/}}  
The confusion described there stems from textbooks emphasizing that p-values and confidence intervals have no probabilistic interpretation {\em and then not offering an alternative interpretation}.
%---only a repeated-sampling-based definition/explanation is given which, to the necessarily statistically-inexperienced textbook reader, must seem wholly unrelated to the practical problem and the data at hand.  
Without an interpretation, there are at least two ways things can go.  Some researchers will craft their own interpretations, but a multitude of different interpretations can only create confusion.  Other researchers will simply accept that there is no meaningful interpretation, making statistical analysis a protocol one blindly follows, i.e., a ``cult of statistical significance'' \citep{cult.stat.significance}.  This confusion and/or blind trust leads to misuse of statistical tools and, perhaps more importantly, encourages researchers to focus on the relatively narrow scientific questions that they believe can be answered with the simple textbook protocols.  Fortunately, this confusion can be overcome, since the frequentist--IM connection offers a simple and mathematically rigorous interpretation of p-values and confidence intervals.  Fisher was correct when he said that p-values and confidence intervals warrant no ``exact probability statements'' about $\Theta$, but that doesn't mean no statements are warranted.  Borrowing Shafer's description of upper probabilities as measures of plausibility, the aforementioned connection immediately implies that p-values can be interpreted as {\em the plausibility of $H_0$, given data $z$} and confidence sets can be interpreted as {\em the set of parameter values that are all individually sufficient plausible, given data $z$}.  This is precisely how p-values and confidence sets are used in practice, and now there's a mathematically rigorous justification for that interpretation.  This is the interpretation of p-values and confidence sets that I teach in my courses, even the introductory-level courses---no need for any technical details about imprecise probability, etc.---and this has been well-received by students.\footnote{For example, a masters student who took my linear models course sent me an email after the course was completed with the following remark after her job interview: ``I wanted to let you know that your discussion of plausibility in ST503 made it easier for me to explain what a p-value is when I was interviewing.  Thank you for the new perspective on the subject!''}  

% Trafimow and Marks: the [p-value] fails to provide the probability of the null hypothesis, which is needed to provide a strong case for rejecting it. 

\subsection{For Bayesians}

In contrast to frequentists, Bayesians are committed to probabilism.  Such a commitment is warranted when genuine prior information is available, but otherwise it's questionable.  Since no prior probability distribution that faithfully represents ignorance, no default-prior Bayesian posterior distribution can be ``correct'' in any sense---``[Bayes's theorem] does not create real probabilities from hypothetical probabilities'' \citep{fraser.copss}.  Moreover, even for the pragmatic Bayesian that isn't concerned about his/her posterior distribution being ``correct'' still must accept the lack of reliability guarantees from the false confidence theorem.  For these and perhaps other reasons, \citet{efron.cd.discuss} wrote:
\begin{quote} \label{quote:efron}
...perhaps the most important unresolved problem in statistical inference is the use of Bayes theorem in the absence of prior information. 
%\citep[][p.~41]{efron.cd.discuss}
\end{quote}
Insisting on probabilism is a constraint on the quality and reliability of uncertainty quantification.  To emphasize this point, the Society of Imprecise Probabilities Theory and Applications has a motto 
%\footnote{See the first sentence on the website \url{https://sipta.org/}.} 
that reads: {\em there is more to uncertainty than probability}.  The IM formulation openly accepts this limitation of probabilism and acknowledges that, while there is no single ``correct'' or fully reliable posterior probability distribution in the absence of prior information, there is a set of posterior probabilities that is justifiably reliable, and that set can be characterized by a possibility measure.  

Understandably, the reader might be uncomfortable with imprecise probability and, for the sake of simplicity, prefer the familiar probabilistic uncertainty quantification despite its shortcomings.  
%For example, at least you know how to summarize a probability distribution, e.g., with moments, quantiles, etc.; see {\color{red} Appendix~??} for a bit about analogous possibilistic summaries. 
But there are many ways to construct probabilities, and focusing solely on ``prior times likelihood'' constructions is again a constraint on the quality of uncertainty quantification.  The new idea in \citet{reimagined} is to {\em approximate} the IM's possibilistic output by a probability distribution.  Below I give an overview.
%, but further details are provided in Appendix~\ref{AA:inner}.  

Recall the inner probabilistic approximation $\prior_z^\star$ which, if it exists, satisfies the condition \eqref{eq:imcd}.  Then such a $\prior_z^\star$ is exactly probability matching in the sense that its credible sets are confidence sets.  Typically, the default-prior Bayes solutions can only achieve this credibility--confidence matching asymptotically.  Moreover, under the usual regularity conditions, $\prior_z^\star$ inherits a probabilistic Bernstein--von Mises theorem from the underlying IM's possibilistic version (Theorem~\ref{thm:bvm}), which implies that the solution $\prior_z^\star$ is asymptotically efficient and, therefore, the exact probability matching isn't due to being overly conservative.  And while the inner probabilistic approximation generally isn't a Bayesian posterior under any prior, there are cases in which a direct Bayesian connection can be made.  In particular, for the so-called invariant statistical models \citep[see, e.g.,][Ch.~6]{eaton1989, schervish1995}, the Bayesian posterior based on the right Haar prior is an inner probabilistic approximation of the possibilistic IM \citep[e.g.,][]{martin.isipta2023, reimagined}.  
%More generally, while the inner probabilistic approximations are not describable in closed-form, is amenable to (approximate) Monte Carlo sampling, as in Section~\ref{SS:computation}.  

To summarize, the limitations of probabilism when prior information is absent, together with the aforementioned benefits of possibilism, warrant abandoning the former in favor of the latter.  But even if one insists on probabilism, the Bayesian-style ``likelihood times prior'' construction itself has limitations: if there was a magical, default prior that solved Efron's problem, then surely that would've been found by now.  So, the solution to that problem will likely come from a completely different perspective, one where the posterior isn't arrived at by application of Bayes's theorem.  Maybe $\prior_z^\star$ is the no-Bayesian-eggs omelet that solves Efron's ``most important unresolved problem''?

\section{Eliminating nuisance parameters}
\label{S:marginal}

%While marginal inference on interest parameters in the presence of nuisance parameters rarely receives careful treatment in textbooks, 
\citet{basu1977} wrote that ``Eliminating nuisance parameters from a model is universally recognized as a major problem of statistics.'' Little has changed since Basu's time---the frequentist impossibility results in, e.g., \citet{gleser.hwang.1987} and \citet{dufour1997}, and the general unreliability of Bayesian inference as discussed above imply that marginal inference is challenging and requires careful treatment.  The possibility-theoretic perspective here offers some new insights, which I discuss below. 

A general operation performed in (imprecise) probabilistic inference is {\em extension}, where uncertainty quantification about one unknown is extended to a related unknown using the uncertainty quantification framework's calculus.  
%\citet[][Ch.~3.1]{walley1991} makes this point clear, even if some of his terminology hasn't been defined here:
%\begin{quote}
%[Extension] is the fundamental concept in our theory of statistical inference... Indeed, natural extension may be seen as the basic constructive step in statistical reasoning; it enables us to create new previsions from old ones. 
%\end{quote}
%Simply put: extension is what inference---statistical or otherwise---is all about.  
In possibility theory, the relevant calculus is optimization, hence this is the operation used to carry out extension.  Following \citet{zadeh1975a, zadeh1978}, the building block of the possibilistic extension principle is optimization-based rule for marginalization: using the present notation and terminology, if $\Theta$ is unknown, with $z$-dependent  quantification of uncertainty given by the possibilistic IM with contour $\pi_z$, and if $\Phi = g(\Theta)$ is a feature of $\Theta$, then the corresponding extension-based marginal IM contour for $\Phi$ is defined as 
\begin{equation}
\label{eq:mpl.ex}
\pi_z^\text{\sc ex}(\phi) = \sup_{\theta: g(\theta)=\phi} \pi_z(\theta), \quad \phi \in g(\TT). 
\end{equation}
%This might seem like the opposite of extension, since the $\Phi$ above is fully determined by $\Theta$.  But equipped with the rule above for marginalization, the process can be reversed to justify a possibility contour for $\Theta$ given only that for a feature $\Phi$ thereof.  Since this proper ``extension'' isn't directly relevant to the present discussion, I defer my explanation of this terminology to Remark~\ref{re:extension} in Appendix~\ref{A:remarks}.  
While the formality leading up to \eqref{eq:mpl.ex} might be unfamiliar, the operation being carried out is one that statisticians use without second thought: to test a composite hypothesis, one can maximize p-values over all the simple hypotheses that it consists of.
%; likewise, given a confidence set for $\Theta$, its image under mapping $g$ is a confidence set for $\Phi=g(\Theta)$.  

Extension-based marginalization is conceptually simple and reliable in the sense that it preserves the IMs strong validity property like in Theorem~\ref{thm:valid}.  That is, for any mapping $g$, if $\pi_z^{\text{\sc ex}, g}$ is the marginal IM contour for $\Phi=g(\Theta)$ as defined in \eqref{eq:mpl.ex}, then 
\[ \sup_{\theta \in \TT} \prob_\theta\{ \pi_Z^{\text{\sc ex}, g}(g(\theta)) \leq \alpha \} \leq \alpha, \quad \alpha \in [0,1]. \]
The simplicity and generality of this strategy are advantageous.  But without any tailoring to the specific problem or feature of interest, one should expect the corresponding extension-based marginal IM to be rather conservative.  
%More efficient marginal IMs require some feature-specific considerations.

When interest is in a specific feature $\Phi=g(\Theta)$, it makes sense to modify the ranking step in the IM construction.  The rationale is that the goal is to rank values $\phi$ of $\Phi$ in terms of their compatibility with the data---how compatible $(z,\theta)$ pairs are is no longer directly relevant.  In line with possibility theory's optimization-centric perspective, a natural strategy is to use a relative profile likelihood:
\[ R^\text{\sc pr}(z,\phi) = \sup_{\theta: g(\theta)=\phi} R(z,\theta), \quad \phi \in g(\TT), \]
where $R$ is the relative likelihood used previously.  Then the validification step as usual completes the construction of a profile-based marginal possibilistic IM for $\Phi$, with contour 
\begin{equation}
\label{eq:mpl.pr}
\pi_z^\text{\sc pr}(\phi) = \sup_{\theta: g(\theta) = \phi} \prob_\theta\{ R^\text{\sc pr}(Z,\phi) \leq R^\text{\sc pr}(z,\phi) \}, \quad \phi \in g(\TT). 
\end{equation}
Note that the outer supremum is needed because, while the relative profile likelihood $R^\text{\sc pr}(Z,\phi)$ only directly depends on $\phi$, its distribution depends on the model parameter, which is not fully determined by $\phi$.  As before, it is easy to show that strong validity is preserved under this profile-based marginal IM construction.  The rationale behind why including an optimization in the ranking step is generally superior in terms of efficiency compared to optimizing post-validification step is subtle and I refer the interested reader to \citet{martin.partial2}
%those details to Remark~\ref{re:complexity.principle} in Appendix~\ref{A:remarks}.  
That profiling tends to be more efficient than extension is easy to see in specific applications; see below.  \citet{imbvm.ext} show that, while both the extension- and profile-based marginal IM constructions enjoy a large-sample possibilistic Bernstein--von Mises theorem, the latter's limiting Gaussian generally has a smaller variance, hence its greater efficiency.  

\begin{figure}[t]
\begin{center}
%\scalebox{0.7}{\includegraphics{hamada_joint_pl}}
\subfigure[Normal example]{\scalebox{0.5}{\includegraphics{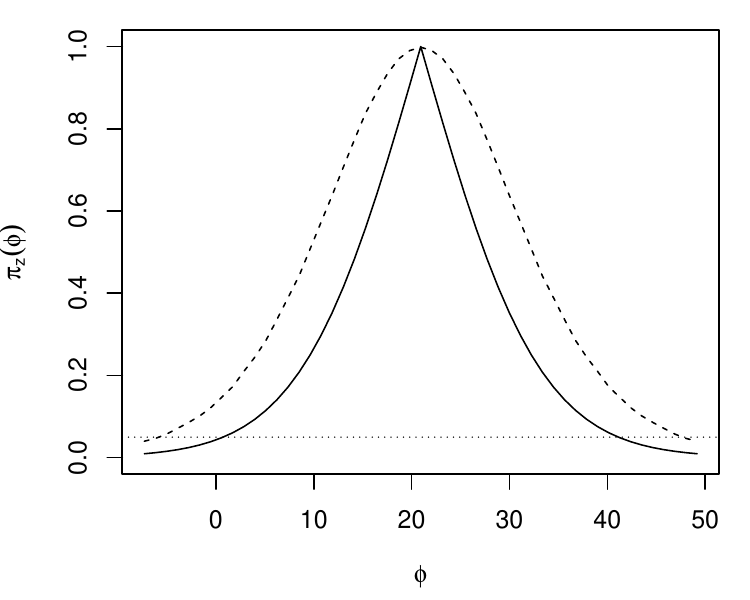}}}
\subfigure[Gamma example]{\scalebox{0.5}{\includegraphics{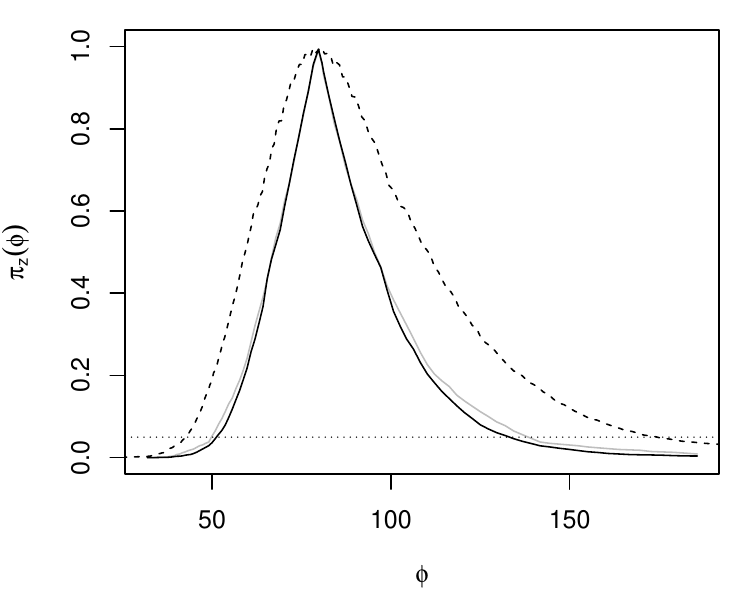}}}
\end{center}
\caption{Plots of the marginal IM contours---extension-based (dashed line) and profile-based (solid line)---for the mean $\Phi$ in the normal and gamma examples. Gray line in Panel~(b) is explained in the text.}
\label{fig:example2}
\end{figure}

\begin{gaussex}[cont]
Here I'll revisit the normal example based on Darwin's data as analyzed in \citet[][Ch.~3]{fisher.doe}.  Here, suppose that the interest parameter is $\Phi = g(\Theta) = \Theta_1$, the mean of the normal distribution.  Since this study involved paired data, the mean $\Phi$ corresponds to the difference between the means of the two marginal populations paired together.  Figure~\ref{fig:example2}(a) shows the extension-based (dashed) and profile-based (solid) marginal IM contours for $\Phi$.  The latter, profile-based contour $\phi \mapsto \pi_z^\text{\sc pr}(\phi)$ exactly corresponds to the p-value for testing $H_0: \Phi=\phi$ based on the two-sided Student-t statistic or, equivalently, the probability-to-possibility transform of the sampling distribution of that Student-t statistic.  Therefore, the value at $\phi=0$ is exactly the usual t-test's p-value which, in this case, is about 0.0495.  Note that the profile-based contour is more tightly concentrated, e.g., the 95\% confidence interval determined by the horizontal line at $\alpha=0.05$ is narrower, hence the profile-based marginal IM is more efficient.   
\end{gaussex}

\begin{gammaex}[cont]
Here I revisit the previous analysis of the data in \citet{hamada.etal.2004}.  Here the focus is inference on the mean $\Phi = \Theta_1\Theta_2$ of the gamma distribution.  Two kinds of marginalization can be carried out---extension-based and profile-based---and both of these are shown Figure~\ref{fig:example2}(b).  As above, the extension-based contour, which is simple to derive from the joint contour shown in Figure~\ref{fig:example1}(b), turns out to be much wider and lacks the efficiency of the profile-based solution.  Figure~\ref{fig:example2}(b) also shows the ``exact'' profile-based marginal IM contour (gray line) based on the brute-force strategy in \eqref{eq:pi.naive}.  I call this ``exact'' because that strategy produces pointwise unbiased estimates of the contour.  I'm showing these two curves here to highlight the accuracy of the sampling-based Monte Carlo strategy of \citet{immc} described briefly in Section~\ref{SS:computation}: the two solid lines in Figure~\ref{fig:example2}(b) are almost indistinguishable.  
\end{gammaex}

Further discussion and illustrations of the profile likelihood-based IM solution are given in Appendix~\ref{A:marginal}.  Despite its advantages over extension-based marginalization, it's important to emphasize that the profiling-based marginalization isn't universal, i.e., there are cases in which profiling is sub-optimal.  Issues arise, as expected, when there are many nuisance parameters, like in the famous examples of \citet{neyman.scott.1948} and of \citet{stein1959}; see \citet[][Sec.~3.6]{martin.partial3}.  More specifically, the profile-based marginal IM is always valid, but the efficiency deteriorates as the number of nuisance parameters increases.  The reason is that maximum likelihood estimators tend to be inconsistent when the number of nuisance parameters diverges; since the profile-based marginal IM contour's peak is at the maximum likelihood estimator, if this is off-target, then a wider contour is needed to cover the relevant range than if the peak was on-target.  The remedy is to replace the relative profile likelihood-based ranking with something else, such as a marginal or conditional likelihood \citep[e.g.,][]{Severini:1994, severini1993, Severini:1998}, but so far this has only been addressed on a case-by-case basis \citep{martin.partial3}.  It deserves to be mentioned a second time that the aforementioned IMs are always valid---unlike Bayes and fiducial, which can be misleading in nuisance parameter problems---so the question is how to properly rank the interest parameter values so that inference is efficient.  

Finally, the ideal cases are those where the relative profile likelihood $R^\text{\sc pr}(Z,\phi)$ is a pivot under $\prob_\theta$ with $g(\theta)=\phi$, since the distributional dependence on the parameters drops out and the computationally challenging supremum on the outside can be omitted.  The relative profile likelihood is a pivot in some cases, and approximately so in many other cases---thanks to Wilks's theorem---but not always.  When that supremum can't be ignored, certain adjustments or approximations may be needed.  Below are two such ``tricks'' to side-step the optimization in \eqref{eq:mpl.pr} that will be used in  Section~\ref{S:beyond} below.
\begin{trick}
One simple and general strategy for eliminating nuisance parameters is conditioning.  A familiar example of this is Fisher's exact test, where the p-value is obtained using the conditional distribution of the test statistic, given the observed value of what is a sufficient statistic under the null.  By definition, the conditional distribution of data given a sufficient statistic doesn't depend on the parameter, hence the nuisance parameter is eliminated.  The present goal isn't getting a p-value to test a hypothesis, but the relevant computations are similar and, therefore, the same strategy can be employed.
\end{trick}
 
\begin{trick}
Strictly speaking, the supremum on the outside of \eqref{eq:mpl.pr} isn't necessary.  Indeed, the practically-out-of-reach contour defined as 
\[ \pi_z^{\text{\sc pr}, \star}(\phi) = \prob_\Theta\{ R^\text{\sc pr}(Z,g(\Theta)) \leq R^\text{\sc pr}(z,\phi) \}, \quad \phi \in g(\TT), \]
which has no supremum, leads
%or, in the case where $\Theta$ can be decomposed as $(\Phi,\Lambda)$, 
%\[ \pi_z^{\text{\sc pr}, \star}(\phi) = \prob_{\phi,\Lambda}\{ R^\text{\sc pr}(Z,\phi) \leq R^\text{\sc pr}(z,\phi) \}, \quad \phi \in g(\TT), \]
%which don't involve suprema, lead 
to a strongly valid marginal IM.  I say that this is 
%these are 
out of reach because, of course, evaluation requires knowledge of $\Theta$ or at least the true nuisance parameter.  But this out-of-reach formulation suggests a simple strategy for constructing a marginal IM that is approximately valid.  That is, for the above case, and for the one where $\theta$ naturally decomposes as $(\phi,\lambda)$, 
%two cases, respectively, 
set
\begin{align*}
\pi_z^{\text{\sc pr}}(\phi) & = \prob_{\hat\theta_z}\{ R^\text{\sc pr}(Z,g(\hat\theta_z)) \leq R^\text{\sc pr}(z,\phi) \} \\
\pi_z^{\text{\sc pr}}(\phi) & = \prob_{\phi,\hat\lambda_z}\{ R^\text{\sc pr}(Z,\phi) \leq R^\text{\sc pr}(z,\phi) \}, 
\end{align*}
where $\hat\theta_z$ and $\hat\lambda_z$ are (maximum likelihood) estimates of the full parameter $\Theta$ and of the nuisance parameter $\Lambda$, respectively.  Strategies of this type have a bootstrap flavor (Section~\ref{S:beyond}) and have been applied successfully by \citet{imcens, immeta} in meta-analysis and survival analysis applications.  
%In fact, that their IM, which is only approximately valid for large samples, was empirically valid even for relatively small samples led these authors to conjecture a sort of higher-order accuracy---beyond the first-order accuracy proved in their papers---for this approximate IM, akin to those established for bootstrap-based methods.  
Improvements to this simple proposal that would be ``closer to valid'' include, for example, optimizing over nuisance parameters in a neighborhood of $\hat\lambda_z$ in the spirit of \citet{berger.boos.1994} and \citet{xiong.optim}.  These improvements will be investigated elsewhere.
\end{trick}

\section{IMs in more general contexts}
\label{S:beyond}

\subsection{Key technical extension}
\label{SS:key}

An apparent limitation of the above proposal is that its emphasis on relative likelihood implicitly assumes a statistical model $\{\prob_\theta: \theta \in \TT\}$ is available.  Machine learning applications, for example, tend to avoid such model assumptions.  A simple but important observation, one that's been applied in various contexts---including Section~\ref{S:marginal} above---is that validification doesn't require the ranking to be based on relative likelihood.  That is, while the posited model's likelihood function determines the ``optimal'' choice of ranking relative to that model (Remark~\ref{re:rellik} in Appendix~\ref{A:remarks}), there may be other factors that suggest a different choice.  Here are a few key instances:
\begin{itemize}
\item If data are coming from multiple sources, e.g., in a meta-analysis or a divide-and-conquer regime \citep{imdc}, and in the form of summary statistics, it may not be possible to compute the full likelihood.
\item More generally, the posited model might not determine a likelihood function for the parameter of interest, e.g., in quantile regression. 
\item If partial/incomplete prior information about $\Theta$ is available---see \citet{martin.partial2} and Section~\ref{S:discuss} below---or if the problem's context suggests some hypotheses are higher priority than others \citep[e.g.,][]{yang.wang.liu.partial, liu.williams.two}, then there's reason to modify the likelihood-based ranking function.
\end{itemize} 
Here I'll briefly describe this simple-but-important technical extension, and then I'll apply it to some relevant problems; see, also, Appendix~\ref{A:general}.  

Relax the assumption that $\Theta$ is the true parameter in a posited statistical model, by saying only that there is a true distribution $\prob$ and that $\Theta \in \TT$ is the true value of the functional $\tau: \model \to \TT$ applied to $\prob$, where $\model$ is some specified subset of $\text{probs}(\ZZ)$, not necessarily a parametric model.  Consider a ranking function $\rho: \ZZ \times \TT \to \RR$ such that large values of $\rho(z,\theta)$ indicate that data $z$ and feature $\theta$ are compatible in some meaningful way.  More detailed example will be given below, but a simple case is where $\rho(z,\theta)$, where $z=(z_1,\ldots,z_n)$ and $z_i = (x_i, y_i)$, represents the negative residual sum of squares corresponding to a fit of some regression function with parameters $\theta$, i.e., 
\[ \rho(z,\theta) = -\sum_{i=1}^n \{ y_i - m_\theta(x_i) \}^2. \]
Note that this setup doesn't assume a parametric model for the data or that the mean structure imposed by $m_\theta$ is ``correct.''  Specification of $\rho$ completes the ranking step of the IM construction, and the validification step proceeds just like in Section~\ref{S:marginal}.  That is, since this is effectively a marginal inference problem in which all but the feature $\tau(\prob)$ of $\prob \in \model$ is a nuisance parameter to be eliminated using the aforementioned recipe:
\begin{equation}
\label{eq:pi.working}
\pi_z(\theta) = \sup_{\prob \in \model: \tau(\prob) = \theta} \prob\{ \rho(Z,\theta) \leq \rho(z,\theta) \}, \quad \theta \in \TT. 
\end{equation}
It is easy to verify that this IM is still strongly valid (relative to $\model$) in the sense that 
\[ \sup_{\prob \in \model} \prob\{ \pi_Z(\tau(\prob)) \leq \alpha \} \leq \alpha, \quad \alpha \in [0,1]. \]
Further details about this general IM formulation are given in Appendix~\ref{A:general}.  It's expected that the Bernstein--von Mises theorem in \citet{imbvm.ext} established for the likelihood-based possibilistic IM can be extended to cover some more general cases like described here, but the details remain to be worked out.  

The challenge, of course, is evaluating the supremum in \eqref{eq:pi.working}, and this boils down to making strategic choices of $\rho$ and/or applying the marginalization tricks in Section~\ref{S:marginal}.  Interesting and practical examples of both cases will be discussed next.

\subsection{Inference on risk minimizers}

Let $Z^n=(Z_1,\ldots,Z_n)$ be iid random variables with distribution $\prob$.  Consider a loss function $(z,\theta) \mapsto \ell_\theta(z)$ that measures how (in)compatible a parameter value $\theta$ is with a data point $z$, where large values correspond to a loss due to incompatibility of $\theta$ and $z$.  Of interest is the minimizer $\Theta$ of the risk, i.e., 
\[ \Theta = \arg\min_\theta r(\theta), \quad \text{where} \quad r(\theta) = \int \ell_\theta(z) \, \prob(dz). \]
The risk---or expected loss---function depends on the true $\prob$ and, of course, so does $\Theta$.  But data $Z^n$ carries relevant information about $\Theta$ and so the goal is to quantify uncertainty about the risk minimizer $\Theta$, given the observed data $Z^n=z^n$. 

Situations like this are common in statistics and machine learning.  Aside from the familiar squared error loss, which is used when the goal is inference on the mean (function), other common examples include classification problems with zero--one loss 
\[ \ell_\theta(z) = 1\{ y = c_\theta(x) \}, \quad z=(x,y), \]
for a posited classification function $c_\theta$ parametrized by $\theta$, and quantile regression problems with the check loss, depending on a specified quantile level $u \in (0,1)$, 
\[ \ell_\theta(z) = \{ y - q_\theta(x) \} \, \bigl[ u - 1\{ y < q_\theta(x)\} \bigr], \quad z=(x,y), \]
where $q_\theta$ is a posited quantile function. 

The standard approach replaces the unknown risk function with its empirical version,
\[ r_{z^n}(\theta) = \frac1n \sum_{i=1}^n \ell_\theta(z_i), \quad \theta \in \TT. \]
Then the empirical risk minimizer, $\hat\theta_{z^n} = \arg\min_\theta r_{z^n}(\theta)$ is the natural estimator, often called an M-estimator in the statistics literature.  Generalizing the strategy based on the relative likelihood above, here I use the empirical risk function to rank the parameter values, thus completing the ranking step: 
\[ \rho(z^n, \theta) = -\{ r_{z^n}(\theta) - r_{z^n}(\hat\theta_{z^n}) \}, \quad \theta \in \TT. \]
Unfortunately, the validification step is highly non-trivial since it involves an optimization over all $\prob$'s having a given value of the risk minimizer.  The recommendation in \citet{cella.martin.imrisk}, which is equivalent to an application of ``Trick~2'' in Section~\ref{S:marginal}, is to replace the true $\prob$ with $\widehat\prob_{z^n}$, the empirical distribution, in the validification step.  That is, they propose the possibilistic IM with contour function 
\[ \pi_{z^n}(\theta) = \widehat\prob_{z^n}\{ \rho(Z^n, \hat\theta_{z^n}) \leq \rho(z^n, \theta)\}, \quad \theta \in \TT. \]
Note that the data-driven plug-in shows up in two places: one is obvious on the outside of the probability statement; the other is more subtle since, inside the probability statement, the ``true'' risk minimizer is replaced by the empirical risk minimizer because now, with respect to $\widehat\prob_{z^n}$, the latter is the ``true'' risk minimizer.  To be clear, in the above display, $Z^n$ is an iid sample from $\widehat\prob_{z^n}$, which is equivalent to sampling with replacement from the observed $z^n$---hence a connection to the bootstrap \citep[e.g.,][]{efron1979, efrontibshirani1993, davison.hinkley.1997}.  So, what's proposed here boils down to the construction of a sort of bootstrap-based marginal possibilistic IM for the risk minimizer.  The cost of plugging in the empirical distribution is that exact validity no longer holds, at least not in general.  But \citet{cella.martin.imrisk} prove an asymptotic validity result and show empirically that the validity property---which is only guaranteed to hold in large samples---typically holds even in rather small samples.  

Admittedly, the solution described above is not fully satisfactory, since validity only holds approximately, as the sample size converges to infinity, instead of exactly in finite samples.  But I think there is plenty of room for improvement here, hence my interest in presenting this less-than-fully-satisfactory solution in the present review.  Indeed, my conjecture is that the asymptotic validity described above is higher-order accurate, i.e., the convergence of $\pi_{Z^n}(\Theta)$ to uniformity is faster than the usual root-$n$ rate.  More generally, I'm certain that variations on the above proposal that are at least ``closer to exactly valid'' are within reach. I hope this review inspires others to contribute their own ideas to solve this important and challenging open problem.

\subsection{Prediction}

To showcase the flexibility of the IM formulation, here I'll consider a different sort of problem, namely, prediction.  Assume that $Z^n=(Z_1,\ldots,Z_n)$ consists of iid observations from a common distribution $\prob$, and that the goal is to predict the next observation $Z_{n+1}$.  The setup here can be generalized in at least two directions: first, all that's required is the $Z$-process is exchangeable; second, if $Z_i=(X_i,Y_i)$, then it's easy to accommodate the case of predicting $Y_{n+1}$ based on a given $X_{n+1}$ and the observed $Z^n$.  It's just for simplicity of presentation that I focus on the simple iid case.  The prediction problem can be viewed as an extreme case of marginal inference, where the entirety of $\prob$ is a nuisance parameter to be eliminated.  A very common method in the statistics and machine learning literature these days is {\em conformal prediction} \citep[e.g.,][]{vovk.shafer.book1, shafer.vovk.2008}.  A close connection between conformal prediction and IMs has already been demonstrated in \citet{imconformal, imconformal.supervised}, and what I present below offers new perspectives. 

Like above, if $z^n$ is the observed data, then the goal is first to rank candidate values $z_{n+1}$ of the next observation in terms of their compatibility or {\em conformity} with $z^n$.  This will be done via a ranking function $\rho: \ZZ^n \times \ZZ \to \RR$ such that large values indicate higher compatibility/conformity.  For example, suppose that $\hat z_{z^n}$ is a point prediction of $Z_{n+1}$ based on data $z^n$; this can be anything, but a reasonable example in this iid case might be $\hat z_{z^n} = n^{-1} \sum_{i=1}^n z_i$, the sample mean.  Then the ranking function would be
\[ \rho(z^n, z_{n+1}) = -d(\hat z_{z^n}, z_{n+1}), \]
where $d \geq 0$ is a distance on $\ZZ$.  The only other constraint on $\rho$ is that it be symmetric in its first argument: that is, $\rho(z^n, z_{n+1})$ doesn't depend on the order of the data $z^n$.  With $\rho$ satisfying the ranking step, the validification step (temporarily) looks like
\[ \pi_{z^n}^\dagger(z_{n+1}) = \sup_\prob \prob\{\rho(Z^n, Z_{n+1}) \leq \rho(z^n, z_{n+1}) \}, \quad z_{n+1} \in \ZZ, \]
where the supremum is over all marginal distributions for the iid $Z$-process or, more generally, over all exchangeable joint distributions.  Note that the right-hand side above involves a probability calculation with respect to the joint distribution of $(Z^n, Z_{n+1})$.  It's easy to verify that this predictive possibilistic IM is strongly valid, i.e., 
\begin{equation}
\label{eq:pred.valid}
\sup_\prob \prob\{ \pi_{Z^n}^\dagger(Z_{n+1}) \leq \alpha \} \leq \alpha, \quad \alpha \in (0,1), \quad \text{all $n \geq 1$}. 
\end{equation}
The problem, again, is that the supremum over $\prob$ in the definition of $\pi_{z^n}^\dagger$ is computationally infeasible.  Fortunately, the structure of the problem allows for a not-so-naive implementation of the IM construction, one that sidesteps the above computational difficulty but without sacrificing on finite-sample strong validity in \eqref{eq:pred.valid}.  

As suggested at the end of Section~\ref{SS:key}, a certain conditioning operation can facilitate the elimination of nuisance parameters.  The key insight is that, under the ``hypothesis'' that $Z_{n+1} = z_{n+1}$, the {\em set of values} $\{z_1,\ldots,z_n,z_{n+1}\}$---or, equivalently, the empirical distribution---is a sufficient statistic.  By conditioning on the sufficient statistic, as in ``Trick~1'' from Section~\ref{S:marginal}, the dependence on the unknown $\prob$ is eliminated, thereby simplifying computations.  The new predictive possibilistic IM contour is given by 
\begin{align*}
\pi_{z^n}(z_{n+1}) & = \sup_\prob \prob\bigl[ \rho(Z^n, Z_{n+1}) \leq \rho(z^n, z_{n+1}) \mid \{z_1,\ldots,z_n,z_{n+1}\} \bigr] \\
& = \frac{1}{(n+1)!} \sum_\sigma 1\{ \rho(z^{\sigma(1:n)}, z_{\sigma(n+1)}) \leq \rho(z^n, z_{n+1})\},
\end{align*} 
where the sum is over all $(n+1)!$ many permutations, $\sigma$, of the integers $1,\ldots,n,n+1$.  Finally, since $\rho$ is symmetric in its first argument, the right-hand side above simplifies:
\[ \pi_{z^n}(z_{n+1}) = \frac{1}{n+1} \sum_{i=1}^{n+1} 1\{ \rho(z_{-i}^{n+1}, z_i) \leq \rho(z^n, z_{n+1})\}, \quad z_{n+1} \in \ZZ, \]
where $z_{-i}^{n+1} = \{z^{n+1}\} \setminus \{z_i\}$.  The reader will surely recognize the right-hand side above as the so-called ``transducer'' or ``p-value'' output produced by the inductive conformal prediction algorithm.  In particular, a result establishing what is equivalent to the prediction validity property in \eqref{eq:pred.valid} can be found in Corollary~2.9 of \citet{vovk.shafer.book1}.  The derivation above, which makes use of conditioning on and sufficiency of the empirical distribution more closely resembles that in \citet{faulkenberry1973} and, more recently, \citet{hoff2023}.  That one can arrive at the very powerful conformal prediction methodology through an IM-driven line of reasoning is quite remarkable.

\section{Conclusion}
\label{S:discuss}

This paper surveys some of the recent developments in possibilistic inferential models (IMs).  Most importantly, IMs offer Bayesian-like fully conditional uncertainty quantification which certain frequentist-like calibration properties that imply, among other things, that tests and confidence procedures derived from the IM's output control frequentist error rates.  Neither of the mainstream approaches to statistical inference are able to achieve both the Bayesian-like and frequentist-like goals, so what distinguishes the IM framework is its reliance on imprecise probability, specifically possibility theory.  Fisher hinted that significance tests and confidence intervals warrant ``no exact probability statements'' but he offered no mathematical explanation for this claim.  By being precise about the role played by imprecision, I'm now able to rectify what \citet{efron1998} playfully referred to as ``Fisher's biggest blunder,'' namely, fiducial inference.  I must also emphasize, again, that embracing imprecision doesn't lower the quality of inference and uncertainty quantification---possibility theory is mathematically and philosophically sound and the imprecision prevents false confidence, keeps us honest.  Moreover, the new possibilistic Bernstein--von Mises theorem ensures that, at least asymptotically, the possibilistic IM solution is efficient.
%; this efficiency is also empirically apparent in finite samples, but a proof of this kind of ``optimality'' presently escapes me.  
Much of the discussion here and in the references cited focuses on uncertainty quantification about the parameters in a statistical model, but Section~\ref{S:beyond} describes the first steps towards pushing IMs beyond this relatively narrow case, making key connections to other fundamental ideas in the literature.  
%So, I expect exciting new developments in the coming years.  

Unfortunately, not every recent development could be included in this review.  Next is a brief list of topics that weren't covered in this review.
First, uncertainty quantification has many uses, and one important application is decision-making.  Following the von Neumann--Morganstern program, the Bayesian framework starts with a loss function, which evaluates the quality of an action for a given parameter value, and then seeks an action that minimizes the expected loss, averaging over parameter values with respect to the posterior distribution.  The fiducial framework proceeds similarly \citep[e.g.,][]{taraldsen.lindqvist.2013}.  Possibilistic IMs using Choquet integration is a new approach to assessing the quality of an action in terms of an upper expected loss, and the corresponding decision-theoretic framework offers certain reliability guarantees that Bayesian and fiducial theories don't \citep{imdec, imdec.isipta}. 

Second, this review focused on quantifying uncertainty about the parameters specific to a given statistical model.  But it's often the case that the model itself is uncertain, and this corresponds to an extreme form of marginal inference---all of the model-specific parameters are nuisance.  Preliminary work on IMs in this context can be found in \citet[][Ch.~10]{imbook} and \citet{martin.model.isipta}.  What's missing in these first IM attempts is a penalty on model complexity.  Bayesians achieve this complexity penalization via prior distributions, and the fiducial efforts control this complexity manually \citep[e.g.,][]{hannig.lee.2009, lai.hannig.lee.2013, williams.hannig.2019, hannig.etal.covariance, hannig.etal.pcr, wei.lee.2023, han.lee.2022, hannig.etal.graphon}.  My view is that penalizing model complexity is motivated by a (prior) belief that the true model is relatively simple, and, while mathematicizing vague, partial beliefs like these is impossible using probability theory, it is easy to do using imprecise probability theory.  So, forthcoming work will show how to treat vague beliefs like ``sparsity'' as incomplete prior information, formulated as an imprecise probability, and incorporate these into the IM construction resulting in provably reliable uncertainty quantification about the model itself.  

Third, I've assumed here that prior information is vacuous.  While this is a standard context in the statistical literature, it's arguably rare that investigators literally know nothing about the quantity they aim to infer.  The trouble is that it's equally rare for the available information to be sufficiently complete that it justifies a particular prior distribution for inclusion in a Bayesian analysis.  The model complexity penalization discussed above is a good example---one might believe a structural assumption like ``sparsity,'' but maybe don't know anything about the structure-specific parameters.  If the only options are to embellish on what's known to create a precise prior distribution, or ignore what's known and assume prior information is vacuous, then the latter is a safe choice.  But the relaxed perspective here on uncertainty quantification offers an alternative path, one that encodes {\em exactly} the available prior information, however vague or incomplete it might be, as an imprecise probability for inclusion in the analysis.  This induces a special type of regularization that allows for efficiency gains while retaining validity.  These details are being developed in the series of working papers by \citet{martin.partial, martin.partial2, martin.partial3}.  

There are far too many open problems to list out here, but below are a few that seem particularly interesting, touching on theory, methods, computation, and applications.
\begin{question}
{\em Which statistical hypotheses are afflicted with false confidence?} There is strong theoretical and empirical support for the claim that false confidence is caused by non-linearity, i.e., it's a consequence of probabilistic marginalization through non-linear functions of the full model parameter.  But a precise characterization of what these hypotheses look like and how severe the affliction will be remains elusive.
\end{question}

\begin{question}
{\em What about IMs constructed based on models learned from training data?}  The IM literature assumes that the model form is given, which is somewhat unrealistic.  It's common in machine learning to use training data to learn about certain aspects of the data-generating process, and then use that partially-trained model for inference and prediction.  In the present context, the ranking and/or the validification step can depend on training data.  What can be said about the reliability of IMs constructed in this way?  
\end{question}

\begin{question}
{\em How to scale up to higher dimensions?}  The validity results presented in Section~\ref{SSS:valid} hold for all sample sizes and all parameter dimensions; the only result that assumes the model is ``low-dimensional'' is Theorem~\ref{thm:bvm} on efficiency.  So, concerning scaling up to higher dimensions, the question boils down to computational and statistical efficiency.  In high-dimensional problems, statistical efficiency is achieved through proper regularization, as discussed above, and work in this direction is in progress.  From a computational perspective, new ideas are needed to leverage cutting edge strategies from both optimization and Monte Carlo integration.  I don't expect that brand new ideas are needed, so a good start would be novel combinations of different ideas.  To be fair to IMs, Bayesians and frequentists have been working on high-dimensional problems for a while, and the relevant computational challenges haven't really been ``solved''---we generally know {\em how} to attempt optimization and sampling in high-dimensional problems, but we generally can't prove that these attempts actually work.  
\end{question}

\begin{question}
{\em What about causal inference, differential privacy, etc?}  Exciting applications nowadays involve causation \citep[e.g.,][]{imbens.rubin.book, pearl.causality.2009}, data privacy considerations \citep[e.g.,][]{garfinkel.dp.book, awan.wang.2025}, etc.  There's no fundamental barrier preventing IMs from contributing in these directions, especially given the extension in Section~\ref{S:beyond} above.  It's just a matter of working out the details. 
\end{question} 

I'll end with some very high-level thoughts about IMs and the potential role that they have to play in artificial intelligence (AI).  Of course, AI is concerned with computational systems having the capability to perform tasks typically associated with human intelligence, such as learning, reasoning, problem-solving, perception, and decision-making.  It's not unreasonable to lump this under the broad umbrella of data-driven uncertainty quantification.  Indeed, some psychologists \citep[e.g.,][]{gigerenzer.murray.book, juslin.etal.2007} model cognitive processes as (intuitive) statistical inference: questions are posed, relevant data are collected, and judgments are made based on data, assumed models, etc.  An AI--IM connection is surely difficult to see at this point, in large part because the IM constructions discussed here are specifically tailored for statistical applications.  But the underlying idea---uncertainty quantification with reliability guarantees---is more general and has broader appeal/applicability.  Indeed, the deep learning models employed in modern AI applications are arguably just ``fancy nonparametric regression models,'' so the IM details discussed in Section~\ref{S:beyond} above and in the supplementary material are clearly relevant.  In any case, just as Shafer's efforts to develop what was later called ``Dempster--Shafer theory'' independent of the probabilistic language and statistical focus in Dempster's early work found real AI applications in the 1980s, I'm optimistic that there's a sufficiently general IM formulation that can meet the ``reliable uncertainty quantification'' needs of modern AI.

\section*{Acknowledgments}

This work is partially supported by the U.S.~National Science Foundation, under grants SES--2051225 and DMS--2412628.  The author sincerely thanks Chuanhai Liu for his support, encouragement, and insights over the last 15+ years.  More specifically, thanks also go out to Leonardo Cella, Emily Hector, Sahil Patel, and James Robertson for helpful suggestions and assistance.

\appendix

\section{Possibility theory primer}
\label{A:possibility}

Possibility measures \citep[e.g.,][]{dubois.prade.book} are among the simplest forms of imprecise probability, closely linked to fuzzy set theory \citep[e.g.,][]{zadeh1978} and Dempster--Shafer theory \citep[e.g.,][]{shafer1976, shafer1987, cuzzolin.book}.  Some applications in statistics are described in \citet{dubois2006}.  Probability and possibility theory differ philosophically and mathematically, but here I'll only briefly discuss the latter; for a discussion of the former, the presentation in \citet{shackle1961} is clear and compelling.   

The mathematical differences between probability and possibility theory can be succinctly summarized as follows: optimization is to possibility theory what integration is to probability theory.  That is, a possibility measure $\uPi$ defined on a space $\TT$ is determined by a function $\pi: \TT \to [0,1]$ with the property that $\sup_{t \in \TT} \pi(t) = 1$.  This function is called the {\em possibility contour} and the supremum-equals-1 property is a normalization condition akin to the integral-equals-1 property of probability densities.  Then the possibility measure is determined by optimizing its contour, i.e., $\uPi(A) = \sup_{t \in A} \pi(t)$, for any $A \subseteq \TT$, just like a probability measure is determined by integrating its density.  

The different calculus has a number of implications.  One is that the aforementioned supremum-equals-1 normalization condition ensures that $\uPi$ is a coherent upper probability \citep[e.g.,][]{cooman.poss1, cooman.aeyels.1999, walley1997} in the spirit of \citet{walley1991} and others.  This means that $\uPi$ determines a non-empty (closed and convex) set of ordinary probabilities that it dominates:
\begin{equation}
\label{eq:credal}
\cred(\uPi) = \{ \prior \in \text{probs}(\TT): \prior(H) \leq \uPi(H) \text{ for all measurable $H$}\}, 
\end{equation}
where $\text{probs}(\TT)$ is the set of probabilities supported on (the Borel $\sigma$-algebra of measurable subsets of) $\TT$. The set $\cred(\uPi)$ is called the {\em credal set} and all (coherent) upper probabilities have an associated credal set.  Aside from being relatively simple, an advantage of possibilistic uncertainty quantification is that the associated credal set has a statistically oriented characterization, which is important and beneficial for statistical applications. 

Finally, if $\uPi$ is interpreted as an upper probability, then there is a dual lower probability, $\lPi$.  This dual is called the necessity measure and is given by $\lPi(A) = 1 - \uPi(A^c)$, $A \subseteq \TT$.  This can be interpreted as a measure of support for or confidence in the truthfulness of an assertion $A$.  
%For example, the $100(1-\alpha)$\% confidence sets $C_\alpha(z)$ defined in \eqref{eq:region} satisfy $\lPi_z\{ C_\alpha(z) \} = 1-\alpha$, i.e., $C_\alpha(z)$ is a set to which the IM's necessity measure $\lPi_z$ assigns $1-\alpha$ support/confidence.  
But since $\lPi$ and $\uPi$ are linked together through the above relationship, it suffices to work with just one and I prefer to work with $\uPi$; this is why the lower probability $\lPi$ hardly appears in this review.  

A possibility-theoretic notion and operation relevant to the IM developments in this review is the {\em probability-to-possibility transform} \citep[e.g.,][]{dubois.etal.2004, hose2022thesis, hose.hanss.2021}.  Consider a random variable $T$ with probability density function $f$ supported on a space $\TT$.  Suppose one wants to approximate the probability distribution $\prob$ that quantifies uncertainty about $T$ by a possibility measure---how would this be done?  The probability-to-possibility transform, i.e., 
\[ \pi(t) = \prob\{ f(T) \leq f(t) \}, \quad t \in \TT, \]
is the best such approximation in a sense to be described.  But, first, notice that the probability-to-possibility transform resembles a p-value calculation.  I hope that this (and the forthcoming explanation of why this kind of construction is ``best'') helps to explain why the possibilistic IM construction above shares some similarities with p-values.  That is, the possibilistic IM construction isn't designed to line up with p-values, that's just what the ``best'' possibilistic approximations happen to look like.  Also notice that what the formula in the above display is doing is assigning possibility values to match certain tail probabilities.  This should help explain why, in Section~\ref{SS:computation} above, I propose to approximate the possibility measure by a probability distribution that agrees with the former in the tails---I'm just inverting the probability-to-possibility transform.  To understand in what sense the possibility measure $\uPi$ with contour $\pi$ as in the above display is the ``best'' approximation of $\prob$, the credal sets are needed.  For $\uPi$ to be a reasonable approximation of $\prob$, it makes sense that the credal set $\cred(\uPi)$ should contain $\prob$.  Moreover, the ``best'' possibilistic approximation should correspond to that having the smallest credal set that contains $\prob$.  It turns out that the possibility measure $\uPi$ determined by the contour in the above display satisfies
\[ \cred(\uPi) \subseteq \cred(\uPi') \quad \text{for all $\uPi'$ with $\prob \in \cred(\uPi')$}. \]
For details, see \citet[][Theorem~3.2]{dubois.etal.2004} and \citet[][Sec.~3.1.4]{hose2022thesis}.  The intuition behind this notion of ``best'' is that, while I'm apparently willing to introduce some degree of imprecision when approximating a probability by a possibility, there's no rationale for introducing more imprecision than is absolutely necessary to capture $\prob$ in the corresponding credal set.  

As an example, Section~\ref{SSS:efficiency} in the main text introduced what \citet{imbvm.ext} called the {\em Gaussian possibility measure}.  This is nothing but the probability-to-possibility transform applied to a given Gaussian probability distribution.  Then the interpretation of the Gaussian possibility is that it's the possibility measure with the smallest credal set that contains the given Gaussian probability distribution.

\section{Connecting old and new IM constructions} 
\label{A:old.vs.new}

The original IM formulation advanced in \citet{imbasics} involves a three-step construction.  The first step, called the {\em association step}, gives an explicit link between the observable data $Z$, uncertain parameter $\Theta$, and unobservable auxiliary variable $U$, which I'll write here simply as 
\[ Z = a(\Theta, U), \]
where $a$ is a known function and $U$ is a random variable with a known distribution, denoted by $\prob^{(U)}$, say.  The interpretation is that the above relationship, together with the knowledge that $U \sim \prob^{(U)}$ describes how data are generated according to the posited model.  Of course, for a given model, there are lots of different ways to simulate from that model, so the association is not unique; more on this below.  This is exactly the starting point for Fisher's fiducial inference, Dempster's generalization, Fraser's structural inference, Hannig's generalized fiducial inference, among others.  Excluding impractically simple problems, all of the aforementioned approaches require some pre-processing before the relationship in the above display can be used for inference, but I won't describe these steps here; suffice it to say that this pre-processing involves making certain transformations, conditioning on features of $U$ that are observed, etc.  The point that I want to make is that all of the above approaches, in one way or another, aim to treat $U$ as a random variable both before and after $Z$ is observed.  Of course, if $U$ is interpreted as the generator of the data, then $U$ is random before $Z$ is observed; but once $Z=z$ is observed, $U$ is more-or-less fixed at a value depending on the unknown $\Theta$.  Fiducial-like approaches argue that there's little difference between being {\em fixed at an unknown point} and being {\em random}, so, for given $Z=z$, they invert the above relationship to write $\Theta$ as a function of the known $z$ and assumed-to-be-random $U$, yielding a probability distribution for $\Theta$.  This is where fiducial-like inference frameworks and the IM frameworks differ.  

The IM framework explicitly acknowledges that there's a fundamental difference in the meaning and form of uncertainty associated with a variable $U$ being {\em fixed-but-unobserved} and being {\em random}.  The point is that, if the above display is a description of the data-generating process, then knowing both $Z$ and $U$ would effectively allow one to solve for $\Theta$.  While $Z=z$ is observable, it is impossible to know the value of $U$ that corresponds to $\Theta$ and $Z=z$.  It is possible, however, to reliably predict\footnote{Arguably a better word for this is ``postdict'' \citep[e.g.,][]{gelman.meng.stern.1996, dempster1971} since the guessing is being done after the value of $U$ is chosen and the value $Z=z$ is observed, but that's an awkward and unfamiliar word so I won't use it here.} where the unobserved $U$ is positioned, but not using ordinary probability: for example, drawing new values $U'$ from $\prob^{(U)}$ wouldn't work because actually hitting the target $U$ is a probability-0 event.  Martin and Liu's second step, the {\em prediction step}, introduces a random set $\S$, taking values in the power set of the $U$-space, that is designed to have positive (if not large) probability of hitting the unobservable target $U$.  Again, the point is that a meaningful/reliable inversion of the relationship in the above display requires that the actual unobserved value of $U$ be pinned down.  This can be achieved with a suitably constructed random set $\S$, but not with random draws from a probability distribution.  This is an important step and offers some insights as to the need for imprecision in statistical inference.  That is, Martin and Liu's {\em prediction step} amounts to quantifying uncertainty about the unobservable $U$, given $Z=z$, using a random set rather than a probability distribution.  And since the distribution of random sets can be described by belief functions, possibility measures, and other forms of imprecise probability, the connection becomes clear.  

Finally, once the association that determines the unobservable $U$ and the random set for predicting said $U$ are specified, the final step in the IM construction, the {\em combination step}, puts the different pieces of information together.  Specifically, for a realization of the random set $\S$, define 
\[ \TT_z(\S) = \bigcup_{u \in \S} \{ \theta: z = a(\theta, u)\}. \]
This is a data-dependent subset of the parameter space and, as a function of the random set $\S$, $\TT_z(\S)$ is also a random set.  Then the distribution of this random set is the vehicle through which inference is carried out.  For example, if $H$ is a relevant hypothesis, then the corresponding upper probability is determined by 
\[ \uPi_z(H) = \prob^{(\S)}\{ \TT_z(\S) \cap H \neq \varnothing \}, \quad H \subseteq \TT, \]
where $\prob^{(\S)}$ is the user-specified distribution of the random set $\S$, designed to meet certain easy-to-satisfy calibration criteria that I won't discuss here; see \citet{imbasics} for details.  If these $\S$-calibration conditions are met, then the IM achieves validity like in \eqref{eq:valid.alt} above, i.e., 
\[ \sup_{\theta \in H} \prob_\theta\{ \uPi_Z(H) \leq \alpha \} \leq \alpha, \quad \alpha \in [0,1], \quad H \subseteq \TT. \]
I've admittedly skipped some of the technical details behind the original IM construction here, but those details can be found in \citet{imbasics, imbook} and the other references cited above.  What's more important here in the context of this review is the new perspective that can be offered: key to achieving the desired goal of reliable (imprecise-)probabilistic inference without prior information is acknowledging that probability theory is inadequate, that there's a fundamental difference between a unknown fixed value and a random value, and that statistical inference involves the former case and hence a more flexible framework for quantifying uncertainty.  All that said, carrying out the above three steps to construct an IM for inference in real applications can be challenging, and that's what motivated the developments surveyed in this review.  

Theorem~4.3 in \citet{imbook} is a key result implying that the random set $\S$ should be {\em nested}, i.e., for any two of its possible realizations, one is a subset of the other.  By ``should be nested'' I mean that for any non-nested random set, there exists a nested random set such that the corresponding IM is at least as efficient; in other words, IM constructions using non-nested random sets are {\em inadmissible}.  Since nested random sets are equivalent to possibility measures, another way to view this IM construction is as a choice to quantify uncertainty about the unobservable auxiliary variable $U$ using a possibility measure, and this perspective was explored in \citet{imposs}.  

Whether one prefers the construction in terms of (nested) random sets or possibility measures, the fact remains that there infinitely many choices of these inputs that result in a valid IM.  This begs an important practical question: {\em how to choose the random set/possibility measure?}  One of the issues leading to the multitude of different and apparently acceptable choices is that there's no unique association between data, parameters, and auxiliary variables.  What Martin and Liu recommended was to reduce the dimension of auxiliary variable as much as possible or, in other words, choose an association that could be expressed in terms of an auxiliary variable of lowest possible dimension.  This dimension-reduction step can make us of various techniques, as described in \citet{imcond}, related to but more general than classical notions like sufficiency, conditioning on ancillary statistics, etc.  Additional opportunities for dimension reduction may be possible when the goal is marginal inference on an interest parameter in the presence of nuisance parameters, as explained in \citet{immarg}.  Unfortunately, and not surprisingly, even after putting in the effort to reduce the dimension of the auxiliary variable as much as possible, there's still no obvious/right/optimal choice of random set.  Early efforts focused on constructing optimal random sets when the inferential goal was to assess a particular hypothesis $H$ about the uncertain $\Theta$.  In certain cases at least, connections could be drawn between this optimal random set and uniformly most powerful tests, but both the former and the latter are restricted to certain kinds of (low-dimensional) hypotheses.  Another downside to these efforts is that the optimal random set for a given $H$ need not be the optimal random set for another $H'$ and, since the goal of the IM framework is to reliably and broadly/coherently quantify uncertainty, having a different IM for each hypothesis kills the breadth/coherence of the overall approach.  

\citet{plausfn, gim} took a different approach.  Rather than specifying an association the fully characterized the data-generating process, he instead formed an association between data, parameters, and a {\em scalar} auxiliary variable.  This extreme form of dimension reduction is possible only because the direct connection to the data-generating process is lost.  But with an appropriate choice of summary---the relative likelihood---all the relevant information in the data about $\Theta$ can be preserved, suggesting that this is not a serious sacrifice.  Indeed, to my knowledge, there are currently no examples showing that the approach based on a (reasonable) generalized association as described above is grossly inferior to that based on the formulation in \citet{imbook}.  To the contrary, there are many examples that are out of reach for the original IM formulation that can be readily handled using the generalized IM formulation \citep[e.g.,][]{imcens}.  

An advantage of the generalized association advanced in \citet{plausfn, gim}is that there's really only one reasonable choice of random set, hence the aforementioned ambiguity is gone.  This creates an opportunity to formulate a theory that parallels the (generalized) fiducial and objective Bayes theories that offer a normative approach to statistical inference, one where the theory itself guides the construction of a solution to new problems.  To me, this was a major advancement.  

The possibilistic IM formulation considered here is technically no different than the generalized IM formulation offered in \citet{plausfn, gim}.  What's new in recent years, and the motivation behind writing this review paper, is the connection to imprecise probability theory in general, and possibility theory in particular.  Without this connection, it's easy to interpret the proposal as a variation on or extension of significance testing, a derivative of classical frequentist hypothesis testing theory.  But the goal is broad, coherent, and reliable uncertainty quantification, and the connection to a formal mathematical framework that fits with these desiderata is critical to the developments.  Imprecise probability and possibility theory also provide guidance on how to solve other problems related to inference, including prediction and formal decision-making, although these topics aren't covered in depth in this review.  That same theory also provides some clues on how one can generalize from the ``vacuous prior'' case in consideration here, and this work is currently under development.  Finally, possibility theory also sheds some light on the optimality of the relative likelihood-based construction advocated for here in the review.  The interested reader is referred to Section~4 in \citet{martin.partial2}.

\section{IM-relevant technical remarks}
\label{A:remarks}

\begin{remark}
\label{re:rellik}
Here I aim to justify the claim made in the main text that the relative likelihood-based IM construction is ``right'' or ``optimal,'' at least in cases when prior information is vacuous and there's no justification for prioritizing one hypothesis about $\Theta$ over another.  For this, I'll appeal to some more general results and principles, related to some of the details presented in the possibility theory review in Appendix~\ref{A:possibility}; I'll assume that the reader is familiar with the notation and terminology from there.  This is based on the arguments in \citet[][Sec.~4.2]{martin.partial2}, which itself closely follows the presentation in \citet[][Sec.~2.3.2.1]{hose2022thesis}.  

What I refer to as the ``ranking function'' in the main text is called a {\em plausibility order} in \citet{hose2022thesis}.  For the general random variable $T \sim \prob$ in $\TT$, write this plausibility order as $\rho: \TT \to \RR$, with larger values of $\rho(t)$ indicating higher plausibility.  For a given $\rho$, there is an optimal possibility measure whose contour is given by $t \mapsto \prob\{ \rho(T) \leq \rho(t)\}$.  Principally, this is based on the {\em Principle of Plausibility} which states that ``what is plausible must be possible.''  Mathematically, this choice of possibility contour is ``optimal'' in the sense that its corresponding credal set is the smallest among all those possibility measures whose credal sets contain $\prob$ and whose contours respect the ordering determined by $\rho$.  That we should desire small credal sets is motivated by the {\em Principle of Expressiveness} or {\em Principle of Maximum Specificity} \citep[e.g.,][]{dubois.prade.1986.set}.  So, when it comes to forming a possibilistic representation, the choice comes down to the ranking function or plausibility order $\rho$.  For this, one can appeal to the {\em Principle of Representation} which says ``what is probable must be plausible,'' which suggests that plausibility orders be determined by probability density/mass functions.  This highlights the probability-to-possibility transform described in the main text and also in Appendix~\ref{A:possibility} as the ``right'' or ``optimal'' choice for approximating a probability distribution by a possibility measure; see, also, Theorem~3.2 in \citet{dubois.etal.2004}.  

The statistical context looks a little different only because we interpret likelihood and density functions differently.  But the likelihood can be viewed as an ``upper density function'' when prior information about $\Theta$ is vacuous, and then the use of the likelihood function for ranking $(Z,\Theta)$ pairs is natural.  One last issue is that the likelihood is ranking $\theta$ values for given $Z$, so there's a sort of ``conditioning'' step that's required.  Motivated by minimal sufficiency and a desire to avoid the ambiguities associated with {\em conditioning} in possibility theory \citep[e.g.,][]{cooman.poss2}, the normalized or relative likelihood is used.  For a more detailed explanation, the reader is referred to \citet[][Sec.~4.2]{martin.partial2}.  
\end{remark}

\begin{remark}
\label{re:uniform}
As stated in the main text, strong validity is indeed stronger than validity; the reason is because strong validity is equivalent to a version of the basic validity property that's {\em uniform-in-hypotheses}.  In words, strong validity implies that it's a small probability event that the IM assigns small upper probability {\em to any true hypothesis}, even those chosen by an adversary and/or depending on data.  In mathematics, 
\[ \sup_{\theta \in \TT} \prob_\theta \bigl\{ \text{$\uPi_Z(H) \leq \alpha$ for some $H$ with $H \ni \theta$} \bigr\} \leq \alpha, \quad \alpha \in [0,1]. \]
This reveals how the IM's uncertainty quantification is more than significance testing.  \citet{cella.martin.probing} argue that this uniformity helps bring to practical use the philosophically appealing notions of {\em probing} laid out in \citet{mayo.book.2018}.  
\end{remark}

\begin{remark}
\label{re:comp}
One very simple strategy for approximating the IM contour function is based on the new possibilistic Bernstein--von Mises theorem presented in Section~\ref{SSS:efficiency}.  This boils down to a simple Wilks-style chi-square approximation:
\[ \pi_z(\theta) \approx 1 - F_d\bigl( -2\log R(z, \theta) \bigr), \quad \theta \in \TT \subseteq \RR^d, \]
where $F_d$ is the $\chisq(d)$ distribution function.  The right-hand side is effectively a closed-form expression, and is often quite accurate, even for moderate or small sample sizes.  Even more sophisticated approximations using Bartlett corrections, etc.~can be employed, as discussed in, e.g., \citet[][Ch.~7]{schweder.hjort.book}.  But large-sample approximations are not fully satisfactory.  Ideally, there'd be a middle-ground between the effective-but-expensive strategy in Equation \eqref{eq:pi.naive} of the main text and the too-simple large-sample approximations above.  That's what the sampling-based strategy in \citet{immc} reviewed in the main text aims to achieve.   
\end{remark}

\section{Conditional, fixed-data properties of IMs}
\label{A:conditional}

Famously and provocatively, \citet[][p.~x]{definetti.vol1} starts his treatise on probability by saying ``Probability does not exist.''  What he means is that probability doesn't exist in any non-subjective sense, that probability is a mental construct that depends on one's subjective assessments and judgments.  He further argues that, while individuals have inherent beliefs, accurately eliciting these degrees of belief is unrealistic.  Therefore, de Finetti focuses on what individuals {\em do} rather than what they {\em think}, and suggests relating probabilities to prices an agent is willing to pay, or accept as payment, for a ticket that pays \$1 if a given uncertain event happens (or a given hypothesis is true) and \$0 otherwise.  
%{\color{red} The point is that probabilities, beliefs, etc aren't real, but behaviors are real...}  
What mathematical properties should an agent's pricing system satisfy?

Dutch book theorems, advanced in \citet{ramsey1926}, \citet{definetti.vol1}, and \citet{savage1972}, state that if, say, my pricing system doesn't satisfy the properties of a (finitely additive) probability, then I can be made a sure-loser, i.e., another agent can find a finite collection of gambles which are all acceptable to me but together will result in me losing money no matter what happens/is true.  
%The converse also holds, i.e., if my pricing system satisfies the axioms of (finitely additive) probability, then I can't be made a sure-loser {\color{blue} (Lehman 1955; Kemeny 1955)}.  
Whether the betting context is real or artificial, quantifying uncertainty via a system susceptible to sure-loss is understood to be irrational.  This makes a compelling argument in favor of probabilism.  Indeed, Bayesians apparently have strong reason to resist my proposed transition from probabilistic to possibilistic uncertainty quantification: taken at superficial face value, since possibility measures aren't finitely-additive, the Dutch book theorems imply that those employing possibilistic reasoning can be made sure losers.  

There is, however, an important condition in these Dutch book theorems that's easy to overlook: they assume that, for each gamble, the agent is willing to {\em buy} and {\em sell} the corresponding \$1 ticket for the stated amount.  Such a condition makes sense under the symmetry resulting from perfect information.  But if, for example, there is a collection of different dice that could be rolled, but I don't know which one specifically, then I should base my buying and selling prices for the gamble that pays \$1 if the die lands on ``6'' on the lower and upper probabilities of landing on ``6'' over the collection of possible dice, respectively. If one accounts for this asymmetry, then the Dutch book theorems can be generalized: an agent whose system consists of bounds on acceptable buying/selling prices, expressed as an imprecise probability satisfying some very mild conditions (see below), is safe from being made a sure-loser \citep[e.g.,][Ch.~2]{walley1991}.  Intuitively, bets based on imprecise probabilities are more cautious than bets based on ordinary/precise probabilities, so if the latter is safe from sure-loss, then the former must be too.  Therefore, imprecise-probabilistic reasoning, despite not being finitely-additive, isn't subject to the criticism implied by a face-value interpretation of those classical Dutch book theorems.  

Returning to the statistical inference context of the present paper, while none of the (imprecise) probabilistic uncertainty quantification perspectives requires subjective prior distributions, all are subjective in de Finetti's sense: whatever distribution for $\Theta$, given $Z=z$, the data analyst chooses is indeed his/her {\em choice} and can't be considered ``objective'' in any meaningful sense.  
%There's an important question of why the data analyst would choose one such probability over another; see Remark~\ref{re:interpretation} in Appendix~\ref{A:remarks}. 
One might then ask how do these different perspectives stack up in terms of de Finetti's subjectivist considerations.  Default-prior Bayes and (generalized) fiducial return ordinary probabilities, so these are safe from sure-loss for any fixed $z$.  In the case of IMs, for fixed $z$, I can construct a pricing system out of my IM output as follows: 
\begin{align*}
\lPi_z(H) & = \text{my supremum buying price for \$$1(\Theta \in H)$, given $z$} \\
\uPi_z(H) & = \text{my infimum selling price for \$$1(\Theta \in H)$, given $z$}.
\end{align*}
Here and elsewhere, $1(E)$ denotes the indicator function of the event $E$. The idea is as follows: if another agent is selling \$$1(\Theta \in H)$ for less than $\lPi_z(H)$, then I'd buy it from him/her; if another agent is willing to buy \$$1(\Theta \in H)$ for more than $\uPi_z(H)$, then I'd sell it to him/her; otherwise, the gamble is ``too risky'' so I'd neither buy nor sell.  Since the possibilistic IM above satisfies the aforementioned ``very mild conditions,'' it follows that, for each fixed $z$, this IM-based pricing system is also safe from sure-loss.  

%In the context of statistical inference, however, the gambling terminology generally isn't used; but see \citet{crane.fpp} and \citet{shafer.betting}.  So, instead of treating the lower and upper probabilities as bounds on prices for bets, I can treat them as my own data-dependent degrees of belief about $\theta$.  That is, small $\uPi_y(A)$ suggests to me that the data $y$ doesn't support the truthfulness of the assertion $A$; similarly, large $\lPi_y(A)$ suggests to me that data $y$ does support the truthfulness of $A$; otherwise, if $\lPi_y(A)$ is small and $\uPi_y(A)$ is large, then data $y$ apparently isn't sufficiently informative to make a definitive judgment about $A$.  Like the gambler who refrains from betting in certain situations, in the latter case it's best/safest if I refrain from making inference on $A$ and consider a less complex assertion and/or collect more informative data.  

I emphasized ``for any fixed $z$'' because there are more extensive notions of coherence that consider potential losses incurred when gambles are made first based on prices determined {\em a priori}, then prices are updated in light of observation $Z=z$, and later debts are settled.  This involves properties of the ``(prior, data) to posterior'' mapping, and such considerations are more interesting when prior information is non-vacuous; these cases are beyond the present review's scope, but see \citet{martin.partial, martin.partial2}.  

What are those ``very mild conditions'' mentioned above?  Analogous to the collection of dice in the motivating example in the main text, define the collection of precise probability distributions compatible with the imprecise probability $(\lPi_z, \uPi_z)$, i.e., 
%\begin{equation}
%\label{eq:general.credal}
\[ \cred(\uPi_z) = \{ \prior_z \in \text{probs}(\TT): \prior_z(H) \leq \uPi_z(H) \text{ for all measurable $H$} \}, \]
%\end{equation}
where $\text{probs}(\TT)$ is the set of all probability measures supported on (the Borel $\sigma$-algebra of) the parameter space $\TT$.  This is called the {\em credal set} \citep[e.g.,][Ch.~5]{levi1980} associated with the imprecise probability pair $(\lPi_z, \uPi_z)$, defined in \eqref{eq:general.credal} and discussed in Section~\ref{SS:computation} in the main text.  The individual elements $\prior_z$ generally depend on $z$ because the upper bound $\uPi_z$ does.  The same credal set can be defined in terms of the lower probability $\lPi_z$, just with the inequality ``$\leq$'' reversed.  Then the aforementioned mild condition is that the credal set be non-empty, that there is at least one precise probability compatible with the given imprecise probability.  The idea is that if there exists a probability $\prior_z$ such that the $(\lPi_z, \uPi_z)$ pricing system is more cautious than $\prior_z$, then the former's safety from sure-loss implies the same of the latter.  For possibility measures in general, which includes the possibilistic IMs in particular, the credal set being non-empty corresponds to a very simple condition concerning the contour, namely, $\sup_\theta \pi_z(\theta)=1$ \citep[e.g.,][Prop.~7.14]{lower.previsions.book}.  As was mentioned in the construction, this holds for the possibilistic IM, hence no sure-loss.  
%Further details concerning the possibilistic IM's credal set are discussed in Section~\ref{S:reimagined}.  

The IM output can also be evaluated with respect to the classical statistical principles, such as the sufficiency, conditionality, and likelihood principles.  The short explanation is that the relative likelihood-based possibilistic IM satisfies the sufficiency principle---what \citet{mayo2014} and others refer to as the {\em weak likelihood principle}---but it does not automatically satisfy the conditionality and likelihood principles.  As \citet{martin.basu} explains, this is not a shortcoming, but I point the reader to that paper for further explanation.

\section{More examples}
\label{A:examples}

\subsection{Binomial regression}
\label{AA:logistic}

\begin{figure}[t]
\begin{center}
\scalebox{0.6}{\includegraphics{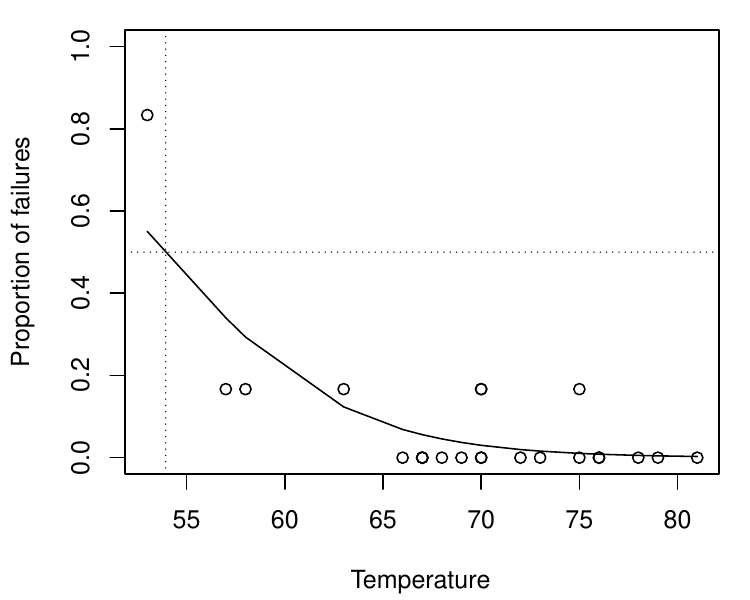}}
\end{center}
\caption{Plot of the o-ring data along with the fitted failure probability curve. The vertical line at 53.94 marks the estimated temperature at which is the o-ring failure probability is 0.5.}
\label{fig:orings.data}
\end{figure}

I'll focus here on a common example involving binomial regression, namely, the o-ring failure experiment associated with the 1986 crash of the Challenger space shuttle.  In this study, 6 tests were conducted at a range of different temperatures ($x$) and the number of o-ring failures ($y$) out of 6 tests were recorded.  The data {\tt orings} are available in, e.g., the {\tt faraway} package \citep{faraway.manual}.  Figure~\ref{fig:orings.data} shows the data and the estimated o-ring failure probability curve obtained by fitting a binomial regression model with logistic link.  Note that the fitted curve is decreasing very fast, suggesting that o-rings are very likely to fail at low temperatures.  This is relevant because the Challenger space shuttle launch took place on a day when the temperature was 31 degrees Fahrenheit.  Here I present some details of an IM analysis of these data.  

\begin{figure}[t]
\begin{center}
\subfigure[Joint contour for $(\Theta_1,\Theta_2)$]{\scalebox{0.55}{\includegraphics{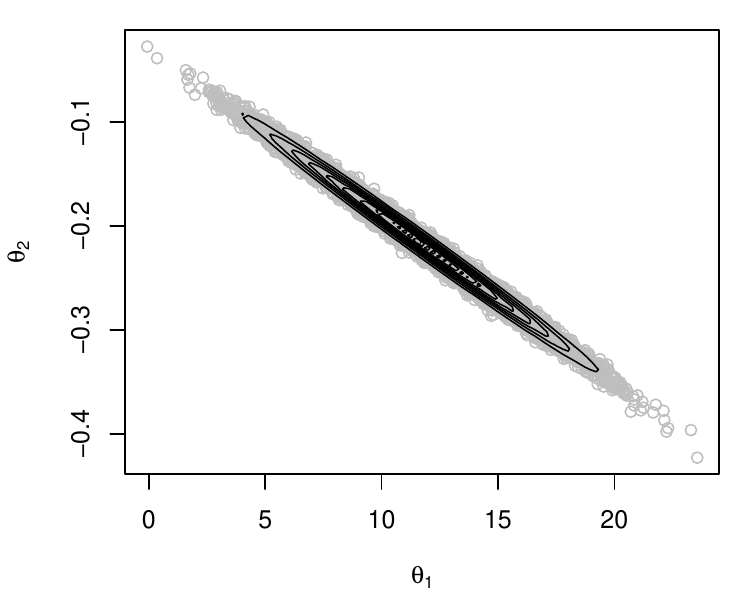}}}
\subfigure[Marginal contour for $\Phi$]{\scalebox{0.55}{\includegraphics{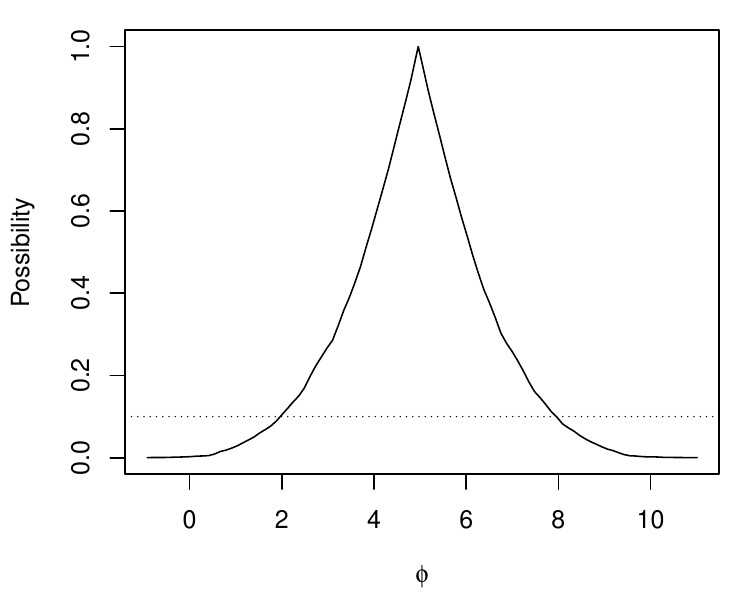}}}
\end{center}
\caption{Panel~(a) shows the samples of $(\Theta_1,\Theta_2)$ from $\prior_z^\star$ and the corresponding approximate possibility contour. Panel~(b) shows the approximate possibility contour for $\Phi$, the log odds of an o-ring failure at 31 degrees.}
\label{fig:orings.out}
\end{figure}

Data $Z$ consists of pairs $Z_i = (X_i, Y_i)$, for $i=1,\ldots,n$, and a conditionally binomial model for $Y_i$, given $X_i$, i.e., 
\[ (Y_i \mid X_i = x_i) \sim \bin\{6, F(\theta_1 + \theta_2 x)\}, \quad i=1,\ldots,n, \]
where where $F(u) = (1 + e^{-u})^{-1}$ is the logistic distribution function.  The corresponding likelihood cannot be maximized in closed-form, but this is easy to do numerically, and the maximum likelihood estimator and the corresponding observed information matrix lead to the asymptotically valid inference reported by standard statistical software.  Exact inference, however, implies a heavier computational burden, since evaluating the exact IM contour over a sufficiently fine grid of $\theta$ values is prohibitively expensive.  Alternatively, the Monte Carlo sampling method presented in Section~\ref{SS:computation} is easy to implement and fast to run.  Figure~\ref{fig:orings.out}(a) shows the 5000 Monte Carlo samples of $(\Theta_1,\Theta_2)$ from $\prior_z^\star$ along with the corresponding approximate possibility contour.  

A practically relevant question concerns the (log) odds of an o-ring failure at temperature 31 degrees, the temperature on the day of the Challenger launch.  Under this binomial model, the log odds of an o-ring failure at 31 degrees is 
\[ \phi = \log \Bigl\{ \frac{F(\theta_1 + \theta_2 \cdot 31)}{1 - F(\theta_1 + \theta_2 \cdot 31)} \Bigr\} = \theta_1 + \theta_2 \cdot 31. \]
For inference on $\Phi = \Theta_1 + \Theta_2 \cdot 31$, the true log odds of failure at 31 degrees, I use the previously-obtained samples $(\Theta_1,\Theta_2)$ from $\prior_z^\star$ to get a corresponding set of samples of $\Phi$.  Figure~\ref{fig:orings.out}(b) plots the corresponding estimate of the possibility contour for $\Phi$.  If the log odds of failure at 31 degrees are roughly between 2 and 8, then that corresponds to a very high probability of o-ring failure, very dangerous.  

\subsection{Categorical data}
\label{AA:mult}

A fundamental problem is inference on the probability vector in a mulinomial model.  I say ``fundamental'' because this is the canonical nonparametric inference problem---one that, e.g., inspired Dempster's original developments \citep{dempster1966} of what's now called Dempster--Shafer theory and still a relevant topic of research, including the discussion paper of  \citep{gong.jasa.mult} published recently in {\em JASA}.  

To set the scene, a random sample of size $n$ is taken from a population consisting of $K$-many categories, labeled $1,2,\ldots,K$.  The model specifies a probability vector $\theta=(\theta_1,\ldots,\theta_K)$ with 
\[ \prob_\theta(\text{observation is of category $k$}) = \theta_k, \quad k=1,\ldots,K. \]
The parameter space $\TT$ is the $K$-dimensional probability simplex, i.e., 
\[ \TT = \{ u \in \RR^K: u_k \geq 0, \, \textstyle\sum_{k=1}^K u_k = 1\}. \]
The observable data $Z=(Z_1,\ldots,Z_K)$ consists of a frequency table listing the number of observations in each of the $K$ categories; note that $\sum_{k=1}^K Z_k = n$.  The true probability vector $\Theta$ is unknown and to be inferred based on observations $Z=z$.  More precisely, let $(Z \mid \Theta=\theta) \sim \mult_K(n, \theta)$.  This determines a likelihood function $L_z(\theta) \propto \prod_{k=1}^K \theta_k^{z_k}$, for $\theta \in \TT$.  Then the relative likelihood is  
\[ R(z, \theta) = \prod_{k=1}^K \Bigl( \frac{n \theta_k}{z_k} \Bigr)^{z_k}, \quad \theta \in \TT. \]
From here it's straightforward to evaluate the possibilistic IM contour function, but expensive to do so over a dense grid of points in the simplex when $K$ is 2 or more.  Fortunately, the computational strategy presented in \citet{immc} and reviewed in Section~\ref{SS:computation} can be readily applied to relieve this computational burden.  In my illustration below, I'll use the basic ellipsoidal approximations as mentioned in Section~\ref{SS:computation} of the main text, but note that this may not be sufficiently accurate when the sample size is small or if there are sparse cells in the observed $z$.  

\begin{figure}[t]
\begin{center}
\scalebox{0.6}{\includegraphics{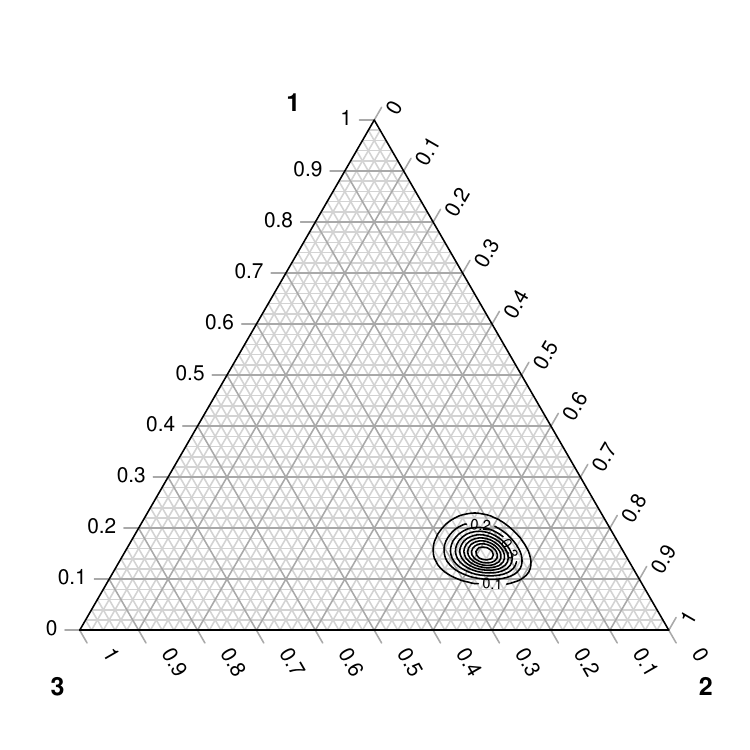}}
\end{center}
\caption{Level sets of the marginal possibilistic IM contour for $\Phi$ in the multinomial illustration described in Appendix~\ref{AA:mult}.}
\label{fig:mult}
\end{figure}

For illustration, I consider a down-scaled version of an example presented in \citet[][Sec.~3.2.2]{agresti2002}.  The original example involves a $3 \times 3$ cross-tabulation of religious beliefs and education levels.  I'm treating this as a multinomial inference problem with $K=9$ categories.  In this form, the original data is 
\[ z^\text{orig} = (178, 138, 108, 570, 648, 442, 138, 252, 252), \]
so that $n^\text{orig} = 2726$.  This is a large sample, making it difficult to visualize the results since the answers are very precise.  So, instead, I down-scale this so the sample proportions are roughly the same but the overall sample size is small; in particular, I use 
\[ z = (10, 8, 6, 34, 38, 26, 8, 15, 15), \]
with corresponding sample size $n=160$.  The full possibilistic IM solution is obtained as mentioned above, resulting in samples of probability vectors from the probabilistic approximation $\prior_z^\star$, which can then be transformed back into a possibility contour.  This relatively high-dimensional solution is difficult to visualize so, for that reason, I'll consider a lower-dimensional feature:
\[ \Phi = (\Theta_1 + \Theta_2 + \Theta_3, \Theta_4 + \Theta_5 + \Theta_6, \Theta_7 + \Theta_8 + \Theta_9). \]
In the original example, this corresponds to the marginal distribution for education level.  The maximum likelihood estimator of $\Phi$ is 
\[ \hat\phi_z = (0.1500, 0.6125, 0.2375), \]
so that category 2 dominates on this margin.  The results of applying the simple probabilistic marginalization that maps the samples of $\Theta$ from $\prior_z^\star$ to samples of $\Phi$, and then transforming those samples of $\Phi$ into a possibility contour are shown in Figure~\ref{fig:mult}.  

\subsection{Variance components}
\label{AA:vc}

Random effect models are common across business, engineering, and science applications.  The simplest such model is often expressed via the following data-generating equation,
\[ Z_{ij} = \lambda + T_i + E_{ij}, \quad i=1,\ldots,m, \quad j=1,\ldots,n, \]
where $Z=(Z_{ij})$ is the observable data and $T_i$ and $E_{ij}$ are mutually independent random effects, with $E_{ij} \iid \nm(0, \phi_1)$ and $T_i \iid \nm(0,\phi_2)$.  Then the three-dimensional parameter that indexes this model is $\theta=(\lambda, \phi_1, \phi_2)$, where $\lambda$ is the overall mean and $(\phi_1,\phi_2)$ is the pair of variance components associated with error/replication and treatment, respectively.  The overall mean is not of primary interest, and it can be easily marginalized out, as shown below.  The focus will be on the two variance components. 

The likelihood function for this model \citep[e.g.,][Eq.~1.4]{tiao.tan.1965} can be expressed as 
\[ L(\lambda, \phi_1, \phi_2) = \phi_1^{-m(n-1)/2} (\phi_1 + n\phi_2)^{-m/2} \exp\Bigl[-\frac12 \Bigl\{ \frac{S_1}{\phi_1} + \frac{S_2 + mn(\bar Z - \lambda)^2}{\phi_1 + n \phi_2} \Bigr\} \Bigr], \]
where $\bar Z=(mn)^{-1} \sum_{j=1}^n \sum_{i=1}^m Z_{ij}$ is the overall sample mean and 
\[ S_1 = \sum_{j=1}^n \sum_{i=1}^m (Z_{ij} \bar Z_{i\cdot})^2 \quad \text{and} \quad \quad S_2 = n \sum_{i=1}^m (\bar Z_{i\cdot}-\bar Z)^2, \]
with $\bar Z_{i\cdot} = n^{-1} \sum_{j=1}^n Z_{ij}$ the group-specific sample mean.  The dependence on $\lambda$ is very simple, so profile likelihood for $(\phi_1,\phi_2)$ can be readily derived:
\[ L^\text{\sc pr}(\phi_1,\phi_2) = \phi_1^{-m(n-1)/2} (\phi_1 + n\phi_2)^{-m/2} \exp\Bigl[-\frac12 \Bigl\{ \frac{S_1}{\phi_1} + \frac{S_2}{\phi_1 + n \phi_2} \Bigr\} \Bigr]. \]
Then the relative profile likelihood $R^\text{\sc pr}(Z,\phi_1,\phi_2)$ is easy to obtain via (numerical) optimization.  It's also clear that the distribution of the relative profile likelihood doesn't depend on the overall mean $\lambda$.  Consequently, the computation 
\[ \pi_z^\text{\sc pr}(\phi_1, \phi_2) = \prob_{\lambda,\phi_1, \phi_2} \{ R^\text{\sc pr}(Z,\phi_1,\phi_2) \leq R^\text{\sc pr}(z, \phi_1, \phi_2)\}, \]
can be carried out by fixing $\lambda=0$, say.  In what follows, I will employ the Bayesian-like sampling strategy using the basic elliptical approximations (applied on the log-variance scale) as described in Section~\ref{SS:computation} of the main text to evaluate the marginal possibilistic IM contour $\pi_z(\phi_1,\phi_2)$ for the variance components.  

\begin{figure}[t]
\begin{center}
\scalebox{0.6}{\includegraphics{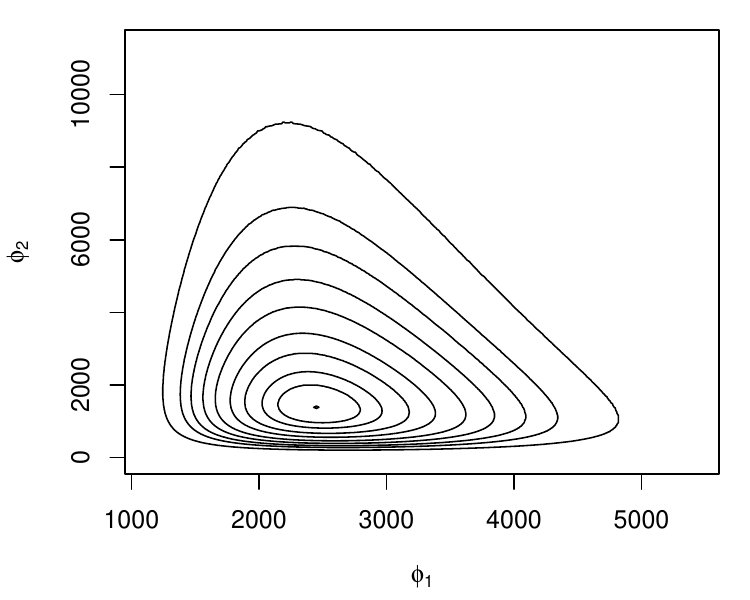}}
\end{center}
\caption{Plots of the marginal possibilistic IM contour for the variance components $(\phi_1, \phi_2)$ as described in Appendix~\ref{AA:vc}.}
\label{fig:vc.pl}
\end{figure}

For illustration, I'll reconsider the real-data example summarized in \citet{tiao.tan.1965}.  The experiment in question was designed to determine if variations in quality of an intermediate product affect the yield of dyestuff prepared from it. In this experiment, there are $m=6$ samples of the intermediate product used and $n=5$ preparations of the dyestuff made from each sample. The observations were measured in grams of standard colors.  Figure~\ref{fig:vc.pl} shows a plot of the marginal IM contour for $(\phi_1,\phi_2)$ based on the data from this experiment.  The point the represents the peak in the IM contour is the maximum likelihood estimator, $(\hat\phi_1, \hat\phi_2) = (2450.7, 1388.6)$.  

%The work that follows is all based on a simulated data set with $m=5$, $n=3$, and true parameters $(\Theta_0, \Theta_1, \Theta_2) = (0, 7, 5)$; the methods that follow depend on the data only through the statistics $S_1 = \sum_i \sum_j (Z_{ij} - \bar Z_i)^2=42.94$ and $S_2 = \sum_i (\bar Z_i - \bar Z)^2 = 37.91$.  A plot of the contour $\pi_z$, as a function of the variance components $(\phi_1,\phi_2)$, is shown in Figure~\ref{fig:vc.pl}; this is based on the new Monte Carlo solution presented in \citep{immc}.  The dot in the center marks the mode of $\pi_z$, the maximum likelihood estimator, $\hat\phi = (4.29, 6.15)$.  

\section{Eliminating nuisance parameters, cont.}
\label{A:marginal}

The goal here is to showcase two different ways in which the (marginal) IM solution is superior to existing solutions when nuisance parameters are present.  

\newtheorem*{ex0}{Example}

\begin{ex0}[Behrens--Fisher]
The classical {\em Behrens--Fisher problem} \citep{fisher1935a, fisher1939} is simple to state: independent samples of size $n_1$ and $n_2$, respectively, are available from two distinct normal populations, $\nm(\Theta_{11}, \Theta_{12}^2)$ and $\nm(\Theta_{21}, \Theta_{22}^2)$, with $\Theta=(\Theta_{11}, \Theta_{12}, \Theta_{21}, \Theta_{22})$ unknown, and the goal is marginal inference on the difference $\Phi = \Theta_{21}-\Theta_{11}$ of the two means. If the two variances are known or if their ratio is known, then the problem is straightforward; but there is no consensus on which solution is ``right'' or ``best'' when the variances are fully unknown.  See \citet{kimcohen1998} for a survey. 

The most common solution modifies the basic Student-t pivot with the degrees of freedom approximations due to \citet{welch1938, welch1947}; this is implemented in R's {\tt t.test} function.  Other standard approaches include the simple-but-conservative solution proposed by \citet{hsu1938} and \citet{scheffe1970}, and the Bayesian solution proposed by \citet{jeffreys1940}, based on the right Haar prior for $\Theta$, which is mathematically equivalent to Fisher's fiducial solution. Ironically, the solution proposed by Jeffreys is different from the Bayesian solution based on the Jeffreys prior, both in its construction and its performance.  A conservative IM solution was advanced in \citet{immarg}, but the approach here, based on the strategy outlined in the main text, is different.  

The approach here is based on validifying the profile relative likelihood.  This was first proposed in \citet{martin.partial3} and an illustration of its use in the Behrens--Fisher problem is considered in Example~5 there.  This is a conceptually simple application of what was summarized in Section~\ref{S:marginal} in the main paper.  The only wrinkle is that the relative profile likelihood doesn't have a closed-form expression in the Behrens--Fisher problem, and its distribution depends on a nuisance parameter; this is precisely what makes the Behrens--Fisher problem challenging.  This increases the computational cost, but I think the validity and efficiency achieved is worth it.  

To showcase the marginal IM's performance, I reproduce a simulation study presented in \citet{reimagined}.  I focus on a rather severely imbalanced case---with $n_1=3$ and $n_2=20$---to ensure that the differences in performance are apparent; this setting was one of the more challenging cases considered in \citet{fraser.wong.sun.2009}.  Otherwise, the simulation settings are standard: $\Theta_{11}=2$, $\Theta_{21}=0$, $\Theta_{12}^2=1$, and $\Theta_{22}^2 = 2$.  I generated 10000 samples under this setup, and Table~\ref{tab:bf} presents the coverage probability and expected length of the various 90\% confidence intervals for $\Phi$.  Notably, only the profile-based marginal IM solution is able to hit the target coverage probability, and, as desired, it's more efficient than the valid-but-conservative Bayes/fiducial solution based on the right Haar prior.  
\end{ex0}

\begin{table}[t]
\begin{center}
\begin{tabular}{ccc}
\hline
Method & Coverage probability & Mean length \\
\hline 
Bayes with Jeffreys prior & 0.861 & 2.31 \\
Bayes with right Haar prior & 0.929 & 3.26 \\
{\tt t.test} with Welch df & 0.881 & 2.79 \\
Profile-based marginal IM & 0.904 & 2.71 \\
\hline 
\end{tabular}
\end{center}
\caption{Coverage probabilities and mean lengths for the 90\% confidence intervals returned by the various methods.  Standard errors for all of the coverage probabilities and mean lengths are roughly 0.003 and 0.01, respectively.}
\label{tab:bf}
\end{table}

\begin{figure}[t]
\begin{center}
\scalebox{0.65}{\includegraphics{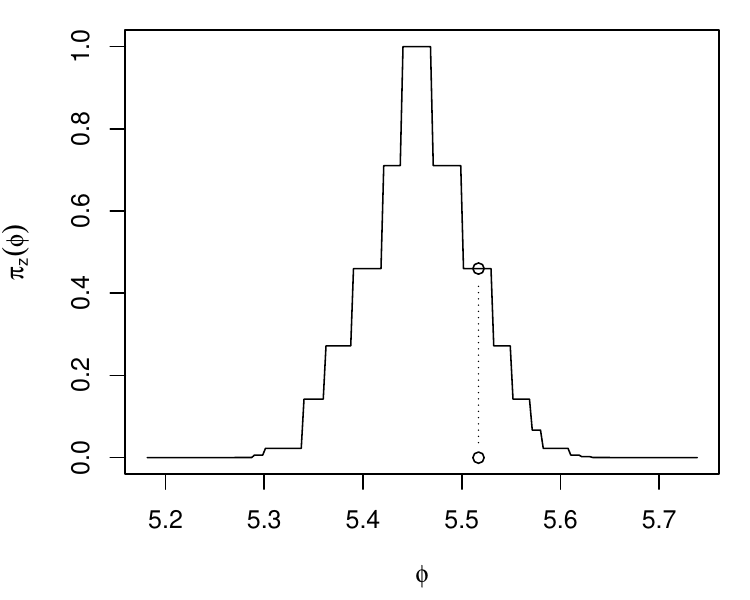}}
\end{center}
\caption{Plot of the marginal IM possibility contour for the median $\Phi$, based on Cavendish's data and the nonparametric relative profile likelihood from \citet{wasserman1990b}.  The vertical line highlights the possibility contour evaluated at the modern ``true'' value, 5.517, of the Earth's density.}
\label{fig:np.med}
\end{figure}

\begin{ex0}[Nonparametric quantile]
This example differs from the others in that it involves a nonparametric problem.  To keep the notation consistent, let $\theta$ be an infinite-dimensional parameter such that $\prob_\theta$, for $\theta \in \TT$, spans the entire class $\text{probs}^\text{\sc ac}(\ZZ)$ of absolutely continuous probability distributions supported on (subsets of) $\ZZ = \RR$; therefore, the true $\Theta$ is in one-to-one correspondence with the underlying true distribution, which has a density with respect to Lebesgue measure on $\RR$.  The data $Z$ consists of $n$ iid observations from distribution $\prob_\Theta$ and the goal is marginal inference on $\Phi = g(\Theta)$, the $q^\text{th}$ quantile of $\prob_\Theta$ for some fixed $q \in (0,1)$.  While there is no likelihood function for $\Phi$ in the usual sense, \citet[][Theorem~5]{wasserman1990b} showed that a relative profile likelihood---or empirical likelihood ratio \citep[e.g.,][]{owen.book}---is well-defined and given by 
\[ R^\text{\sc pr}(z, \phi) = \Bigl( \frac{q}{v} \Bigr)^v \Bigl( \frac{1-q}{n-v} \Bigr)^{n-v}, \]
where 
\[ v = v_q(z,\phi) = \begin{cases} \sum_{i=1}^n 1(Z_i \leq \phi), & \phi \leq \hat\phi_z \\ \sum_{i=1}^n 1(Z_i < \phi), & \phi > \hat\phi_z, \end{cases} \]
and $\hat\phi_z$ is the $q^\text{th}$ sample quantile.  Following \citet[][Example~6]{plausfn}, I propose to construct a marginal possibilistic IM contour $\pi_z^\text{\sc pr}(\phi)$ exactly as in \eqref{eq:mpl.pr}.  For what it's worth, a very similar nonparametric IM construction can be made based on the so-called {\em marginal likelihood} derived in \citet{kalbfleisch1978}, but I didn't find any significant differences between these two solutions.  

For illustration, consider the $n=29$ measurements of the density of the Earth relative to water, taken by Cavendish back in 1798 \citep[][Table~8]{stigler1977}.  These experimental results were used to indirectly calculate the gravitational constant in Newton's formula for the force of gravitational attraction between two bodies.  Here I'll treat the true density of the Earth as the median $\Phi$ of the population of such density measurements.  As above, I construct a marginal IM contour for $\Phi$ based on Cavendish's data, and the result is plotted in Figure~\ref{fig:np.med}.  The modern ``true'' value of the Earth's density is 5.517 and it's clear from the plot that this true value is easily contained in the marginal IM's 95\% confidence interval for $\Phi$.  
\end{ex0}

Related to the second example above, some general comments can be made concerning how the IM's uncertainty quantification framework differs from that of Bayes and fiducial.  In particular, the quantile example reveals that the IM framework can, at least in certain cases, make direct inference about the target quantity of interest.  From the remark quoted in Section~\ref{SS:frequentists} of the main text, \citet{wasserman.quote} continues:
\begin{quote}
%The idea that statistical problems do not have to be solved as one coherent whole is anathema to Bayesians but is liberating for frequentists. 
To estimate a quantile, an honest Bayesian needs to put a prior on the space of all distributions and then find the marginal posterior. The frequentist need not care about the rest of the distribution and can focus on much simpler tasks.
\end{quote}
Of course, Bayesians \citep[e.g.,][]{ferguson1973, ghosal.vaart.book} and (generalized) fiducialists \citep[e.g.,][]{hannig.gfid.survival} can solve nonparametric problems, but they do so by first making inference on the full model parameter---in this case, the unknown distribution function---and then marginalizing to the quantity of interest.  As Wasserman mentions, frequentists can address the marginal inference problem directly, and apparently so too can IM users, at least in this quantile example.  The point is that one doesn't have to give up on data-dependent, fully conditional uncertainty quantification in exchange for a more appearling direct attack on the relevant unknowns.

%\section{ReIMagining details}

%\subsection{IM characterization of frequentist inference}
\section{IM characterization of frequentism}
\label{AA:char}

The possibilistic IM formulation presented in the main text is intended to be normative; that is, I'm offering a rather general construction that is provably valid and at least approximately efficient in a wide range of applications.  I started with a basic likelihood-based construction in Section~\ref{S:details} and that was generalized a bit in Section~\ref{S:beyond}.  It's easy, however, to confuse normativity for rigidity, to think that ``possibilistic IMs'' are only those obtained through those constructions.  From that perspective, one would surely believe that frequentist inference is more flexible than IMs, since there are effectively no rules dictating how the latter is carried out.  But remember that a possibilistic IM is simply a mapping from data $z$ (perhaps depending on a posited model, available prior information, or other things) to an imprecise-probabilistic output whose upper probability has the mathematical form of a possibility measure.  There are, of course, lots of different ways this can be achieved, and fewer---but still many---that would satisfy the validity requirement.  Some of these will be specifically good at answering very specific questions, but lousy at answering other questions.  The point of focusing on the particular likelihood-based construction here and elsewhere in the literature is to offer a solution that's generally good across all relevant questions.  What might be surprising to the reader is that there are exactly as many valid, possibilistic IMs as there are frequentist solutions, as I show formally below.  This justifies the claim made in the main text, namely, that there are no genuine frequentist solutions---including classical those in textbooks as well as those that haven't even been developed yet---that are out of reach by the possibilistic IM framework.  What's notable about the result below is that, while the given frequentist procedure might be a test of a specific null hypothesis or a confidence set for a feature parameter $\Phi$, the resulting possibilistic IM offers full uncertainty quantification about $\Theta$.  The possibilistic IM constructed in the proof below to agree with the given frequentist procedure may not answer every question about $\Theta$ efficiently, but that's not relevant; after all, the frequentist procedure it starts with also can't answer every question efficiently.  The point is simply that whatever frequentist priorities one might have, there's no loss of generality or efficiency focusing on possibilistic IMs.  

Versions of the characterization result presented below have been given in \citet{impval} and \citet{imchar}, with increasing generality.  The result presented below is embarrassingly simple, and requiring only a basic ``nestedness'' condition on the frequentist procedure in consideration.  This amounts to no practical restriction on the characterization result because non-nested confidence sets or tests with non-nested rejection regions are difficult to interpret and rarely if ever used.  

\begin{prop0}
%\label{prop:char}
Consider a model $\{\prob_\theta: \theta \in \TT\}$ for observable data $Z$ with uncertain true parameter value $\Theta \in \TT$.  
\begin{itemize}
\item[(a)] Let $H_0 \subset \TT$ be a relevant null hypothesis and consider a test with a nested family of rejection regions $\{\mathcal{R}_\alpha: \alpha \in [0,1]\}$, i.e., $\mathcal{R}_\alpha \subset \mathcal{R}_{\alpha'}$ for all $\alpha < \alpha'$, such that 
\[ \sup_{\theta \in H_0} \prob_\theta(Z \in \mathcal{R}_\alpha) \leq \alpha, \quad \alpha \in [0,1]. \]
Then there exists a strongly valid possibilistic IM for $\Theta$, with contour $\pi_z$ and corresponding upper probability $\uPi_z$, such that 
\[ z \in \mathcal{R}_\alpha \iff \uPi_z(H_0) < \alpha. \]
\item[(b)] Let $\{C_\alpha^f(z): \alpha \in [0,1]\}$ be a family of confidence sets for $\Phi=f(\Theta)$, nested in the sense that  $C_\alpha^f(z) \supset C_{\alpha'}^f(z)$ for all $\alpha < \alpha'$ and non-empty in the sense that $\bigcap_\alpha C_\alpha^f(z) \neq \varnothing$, with the property that 
\[ \sup_\theta \prob_\theta\{ C_\alpha^f(Z) \not\ni f(\theta) \} \leq \alpha, \quad \alpha \in [0,1]. \]
Then there exists a strongly valid possibilistic IM for $\Theta$, with contour $\pi_z$, such that 
\[ f(\theta) \in C_\alpha^f(z) \iff \pi_z(\theta) > \alpha, \quad \alpha \in [0,1]. \]
\end{itemize} 
\end{prop0}

\begin{proof}
Starting with Part~(a), define the possibility contour $\pi_z$ as 
\[ \pi_z(\theta) = \begin{cases} \sup\{\beta: z \not\in \mathcal{R}_\beta\} & \theta \in H_0 \\ 1 & \theta \not\in H_0. \end{cases} \]
Clearly $\pi_z$ satisfies the normalization condition required to call it a possibility contour.  Then define the corresponding possibility measure $\uPi_z$ via optimization as usual:
\[ \uPi_z(H) = \begin{cases} \sup\{\beta: z \not\in \mathcal{R}_\beta\} & H \subseteq H_0 \\ 1 & H \not\subseteq H_0. \end{cases} \]
Then it's easy to see that 
\[ \uPi_z(H_0) < \alpha \iff \sup\{\beta: z \not\in \mathcal{R}_\beta\} < \alpha \iff z \in \mathcal{R}_\alpha. \]
It's similarly easy to see that 
\[ \prob_\theta\{ \pi_Z(\theta) < \alpha \} = \begin{cases} \prob_\theta(Z \in \mathcal{R}_\alpha) & \theta \in H_0 \\ 0 & \theta \not\in H_0. \end{cases} \]
From this and the assumed Type~I error rate control property of the given frequentist test, the strong validity claim follows immediately.  

For Part~(b), start with defining a possibility contour $\pi_z^f$ for $\Phi$ as follows: 
\[ \pi_z^f(\phi) = \sup\{ \beta: C_\beta^f(z) \ni \phi\}, \quad \phi \in f(\TT). \]
The non-emptiness of the given confidence sets implies that there exists a point at which $\pi_z^f$ equals 1, hence the normalization condition is satisfied and it's fair to call this a ``possibility contour.''  Now extend this to a possibility contour for $\Theta$ in the natural way:
\[ \pi_z(\theta) := \pi_z^f\bigl( f(\theta) \bigr), \quad \theta \in \TT. \]
Then it's clear that $f(\theta) \in C_\alpha^f(z)$ if and only if $\pi_z(\theta) > \alpha$ and, moreover, strong validity holds because 
\[ \prob_\theta\{ \pi_Z(\theta) \leq \alpha \} = \prob_\theta\{ C_\alpha^f(Z) \not\ni f(\theta) \} \leq \alpha, \]
where the right-most inequality holds by the assumed coverage properties of the $C_\alpha^f$'s. 
\end{proof}

%\subsection{Bayes and probabilistic approximations of IMs}
%\label{AA:inner}

\section{General IM construction and properties}
\label{A:general}

The main text focused primarily on IM constructions for a finite-dimensional parameter indexing the posited statistical model.  The reason for this focus is purely practical: most real-world applications are of this type.  But these are not the only cases that the IM framework can handle, and the later sections of the paper briefly mentioned various extensions.  Since this review paper is intended both to showcase the developments that have been made so far and to highlight various opportunities for the future of IMs, I think it's worth sketching some of these details, even if it's a bit beyond the paper's scope.  

The general setup in Section~\ref{S:beyond} of the main paper suggested a possibilistic IM contour $\pi_z(\theta)$ for the relevant feature $\Theta = \tau(\prob)$ of the underlying distribution $\prob \in \model$ as 
\[ \pi_z(\theta) = \sup_{\prob \in \model: \tau(\prob) = \theta} \prob\{ \rho(Z,\theta) \leq \rho(z,\theta) \}, \quad \theta \in \TT. \]
Of course, the corresponding possibility measure or upper probability is defined via optimization just like before: $\uPi_z(H) = \sup_{\theta \in H} \pi_z(\theta)$, for any $H \subseteq \TT$.  Accurate and efficient computation of the IM contour in this very general case is an important open question, so I don't have anything to say about this point here.  My focus in this present discussion is on what kind of properties the ranking function $\rho$ must satisfy in order for the IM constructed above to be valid and efficient.  

One basic requirement is that $\pi_z$ as constructed above should be such that $\sup_\theta \pi_z(\theta) = 1$ for each $z$.  A sufficient condition is that the maximum $\max_\theta \rho(z,\theta)$ be independent of $z$.  This can be achieved, for example, if $\rho$ is such that $\theta \rho(z,\theta)$ is a possibility contour for each $z$---that is, if $\sup_\theta \rho(z,\theta) = 1$ for each $z$.  The relative likelihood and its variations satisfy this condition, but so do many other functions.  To simplify the discussion that follows, I'll assume that $\rho$ takes values between 0 and 1, again, like the relative likelihood and its variants.  

Validity, too, only requires minimal assumptions on $\rho$.  In particular, the distribution of the (scalar) random variable $\rho(Z,\theta)$ should not have an atom at 0 under any distribution $\prob$ with $\tau(\prob) = \theta$.  The reason is that this will imply that $\pi_z(\theta)$ is bounded away from 0 and, therefore, the random variable $\pi_Z(\theta)$ can't be stochastically no smaller than $\unif(0,1)$.  Any practically reasonable choice of ranking function would achieve this since, otherwise, there'd be positive probability of observing data that's fully incompatible (relative to $\rho$) with the true parameter, which clearly makes $\rho$ a poor measure of compatibility.  In any case, mathematically, no atom at 0 in the sense above implies that the distribution function $G_\prob$ of $\rho(Z,\tau(\prob))$ under $Z \sim \prob$ satisfies $\sup_{\prob \in \model} G_\prob(0) = 0$. 

Concerning efficiency, let me turn to large-sample considerations.  For this, write $Z^n$ for an iid sample of size $n$ from unknown distribution $\prob$.  A first basic question is what does it take for the possibilistic IM to be {\em consistent}, i.e., 
\[ \uPi_{Z^n}(H) \to 0 \quad \text{in $\prob$-probability as $n \to \infty$} \quad \text{for each $H \not\ni \tau(\prob)$}. \]
This is analogous to Bayesian (and fiducial) consistency as explored in, e.g., \citet{bsw1999}, \citet{ggv2000}, and \citet{walker2004a}.  To achieve this, the ranking function must be able to distinguish ``good'' and ``bad'' values of the interest parameter $\theta$ as efficiently as possible.  In light of the above constraints, and the exponential growth rate that can be achieved in many problems \citep[e.g.,][Ch.~11]{cover.thomas.book}, there's effectively no choice but to take $\rho$ of the form 
\[ \rho(z^n, \theta) = \frac{\prod_{i=1}^n f(z_i, \theta)}{\sup_\vartheta \prod_{i=1}^n f(z_i, \vartheta)} = \exp\Bigl[ \sum_{i=1}^n \{ \log f(z_i, \theta) - \log f(z_i, \hat\theta_{z^n}) \} \Bigr], \]
where $\hat\theta_{z^n} = \arg\max_\vartheta \prod_{i=1}^n f(z_i, \vartheta)$, for suitable function $f$.  So, either directly or indirectly, explicitly or implicitly, $\Theta=\tau(\prob)$ corresponds to the maximizer of the function $\theta \mapsto \E \log f(Z_1,\theta)$.  This is the case with the relative likelihood as discussed in the main text, since the true parameter is the maximizer of the expected log density function.  

Let $\Theta = \tau(\prob)$ be the true value of the interest parameter and take $\theta \neq \Theta$; then $\theta$ is a ``wrong'' value of the parameter and the IM is expected to rule this value out as the data grows more and more informative, i.e., $\pi_{Z^n}(\theta) \to 0$ in $\prob$-probability as $n \to \infty$; a result of this type was established in \citet[][Theorem~3]{plausfn}.  Towards this, since $\hat\theta_{Z^n}$ is, by definition, a maximizer, it follows that 
\[ \rho(Z^n, \theta) \leq \exp\Bigl[ \sum_{i=1}^n \{ \log f(Z_i, \theta) - \log f(Z_i, \Theta) \} \Bigr]. \]
Under mild integrability conditions, the law of large numbers implies 
\[ \frac1n \sum_{i=1}^n \log\frac{f(Z_i, \theta)}{f(Z_i, \Theta)} \to \E \log \frac{f(Z_1,\theta)}{f(Z_1, \Theta)}, \quad \text{in $\prob$-probability}. \]
Note that the right-hand side above is negative by definition.  Then, by monotonicity of the distribution function, 
\[ \pi_{Z^n}(\theta) = \sup_{\prob \in \model: \tau(\prob) = \theta} G_\prob\{ \rho(Z^n, \theta)\} \leq \sup_{\prob \in \model: \tau(\prob) = \theta} G_\prob\bigl\{ e^{\sum_{i=1}^n \{ \log f(Z_i, \theta) - \log f(Z_i, \Theta) \}} \bigr\}. \]
Clearly, the argument to the distribution function behaves like $e^{-n \gamma}$ for some $\gamma > 0$, hence is converging to 0 very rapidly.  Since there's no atom at 0, the entire right-hand side is vanishing in $\prob$-probability.  Therefore, $\pi_{Z^n}(\theta) \to 0$ in $\prob$-probability for any $\theta \neq \Theta$.  Extension to the proper consistency property as stated above involves establishing that the aforementioned property holds appropriately uniformly in $\theta$.  I won't get into these details here but surely this is doable and would make an interesting contribution.  

Beyond the basic consistency result discussed above, it would be desirable to have concentration rate and generalizations of the Bernstein--von Mises property due to \citet{imbvm.ext} reviewed in the main text.  These details are beyond the scope of high-level discussion here, but these too are deserving of further exploration.  

Finally, as I mentioned briefly in Section~\ref{S:discuss} of the main text, there are even more general IM constructions that can incorporate any available, perhaps incomplete prior information about the relevant unknowns.  These details are too new and, therefore, beyond the scope of this review paper, but they can be found in the working papers \citet{martin.partial, martin.partial2, martin.partial3}.  This incorporation of prior information is done in a novel, non-Bayesian way that allows for efficiency gains akin to those obtained by more familiar forms of regularization but while retaining a generalized notion of validity that's meaningful relative to the (incomplete) prior knowledge available.  A conjecture is that this non-Bayesian style of prior-to-posterior updating is more efficient in the sense that it can avoid the undesirable dilation phenomenon that plagues generalized Bayes solutions \citep[e.g.,][]{kyburg1987, gong.meng.update, walley1991}.

\bibliographystyle{apalike}
\bibliography{/Users/rgmarti3/Dropbox/Research/mybib.bib}

\end{document}